\documentclass[amscd,amssymb,verbatim,12pt]{amsart}
  \usepackage{epsfig}
\newif\ifAMS
\IfFileExists{amssymb.sty}
  {\AMStrue\usepackage{amssymb}}
  {\usepackage{latexsym}}
\theoremstyle{plain}

\newtheorem{Thm}{Theorem}
\newtheorem{Cor}[Thm]{Corollary}
\newtheorem{Lem}[Thm]{Lemma}

\theoremstyle{definition}
\newtheorem{Def}{Definition}

\newtheorem{Se}{Setting}

\theoremstyle{remark}

\newcommand{\interior}{^{ \kern-5pt ^\circ}}
\newcommand {\bd}{\partial}
\newcommand {\G}{\Gamma}
\newcommand {\iy}{\infty}

\newcommand {\N}{{\mathbb N}}
\newcommand {\R}{{\mathbb R}}
\newcommand {\Z}{{\mathbb Z}}

\newcommand {\E}{{\mathbb E}}
\newcommand {\HH}{{\mathbb H}}

\newcommand {\Min}{\text{Min}\,}
\newcommand {\nb}{\text{Nbh}}

\newcommand {\cP}{{\mathcal  P}}
\newcommand {\cR}{{\mathcal  R}}

\begin{document}
\title{Boundaries and JSJ decompositions of $CAT(0)$-groups}

\author
{Panos Papasoglu and Eric Swenson }


\email [Panos Papasoglu]{panos@math.uoa.gr} \email [Eric
Swenson]{eric@math.byu.edu}

\address
[Panos Papasoglu] {Mathematics Department, University of Athens,
Athens 157 84, Greece }
\address
[Eric Swenson] {Mathematics Department, Brigham Young University,
Provo UT 84602}
\thanks {This work is co-funded by European Social Fund and National
Resources (EPEAEK II) PYTHAGORAS}
\begin{abstract}

Let $G$ be a one-ended group acting discretely and co-compactly on
a CAT(0) space $X$. We show that $\bd X$ has no cut points and
that one can detect splittings of $G$ over two-ended groups and
recover its JSJ decomposition from $\bd X$.

We show that any discrete action of a group $G$ on a CAT(0) space
$X$ satisfies a convergence type property. This is used in the
proof of the results above but it is also of independent interest.
In particular, if $G$ acts
co-compactly on $X$, then one obtains as a Corollary that if the Tits
diameter of $\bd X$ is bigger than $\frac {3\pi} 2$ then it is infinite and $G$ contains a free subgroup of rank 2.

\end{abstract}
\maketitle
\section{Introduction}

The purpose of this paper is to generalize results about
boundaries and splittings of hyperbolic groups in the case of
$CAT(0)$ groups. Before stating our results we summarize what is
known in the hyperbolic case.

Bestvina and Mess (\cite{BES-MES}) showed that if the boundary of
a one-ended hyperbolic group does not have a cut point then it is
locally connected. They asked whether the boundary of a one-ended
hyperbolic group can contain a cut point. The negative answer to
this is mainly due to Bowditch (\cite {BOW}, \cite {BOW1}, \cite
{BOW5}) with contributions by Swarup and Levitt ( \cite{Sw},
\cite{Le}). Bowditch \cite {BOW} further showed that the boundary
of a one-ended hyperbolic group has local cut points if and only
if either the group splits over a 2-ended group or it is a
hyperbolic triangle group. He deduced from this a canonical JSJ
decomposition for hyperbolic groups (compare \cite{Se}).

We remark that in the CAT(0) case, boundaries are not necessarily
locally connected (e.g. if $G=F_2\times \Bbb Z$, where $F_2$ is
the free group of rank 2, the boundary is a suspension of a Cantor
set, hence it is not locally connected). However the question
whether the boundaries have cut points also makes sense in this case.
Indeed the second author showed that boundaries have no cut
points under the assumption that the group does not contain an
infinite torsion subgroup (\cite {SWE}).

To state our results we recall some terminology. An action of a
group $G$ on a space $X$ is called {\em proper} if for every
compact $K \subset X$, the set $\{g \in G : g(K) \cap K \neq
\emptyset \}$ is finite. An action of a group $G$ on a space $X$
is called {\em co-compact} if the quotient space of $X$ by the action of
$G$,  $X /G$ is compact.   An action of a group $G$ on a metric
space $X$ is called {\em geometric} if $G$ acts properly, co-compactly
by isometries on $X$.  It follows that $G$ is quasi-isometric to
$X$, that is they have the same coarse geometry. If $Z$ is a
compact connected metric space we say that $c$ is a {\em cut point} of
$Z$ if $Z-c$ is not connected. We
say that a pair of points $\{a, b \}$ is a {\em cut pair} of $Z$  if  $Z- \{a,b\}$ is not connected.

We show the following.
\begin{Thm} Let $G$ be a one-ended group acting
geometrically on the CAT(0) space $X$. Then $\bd X$ has no cut
points.
\end{Thm}
We note that  Croke-Kleiner (\cite{C-K}) showed it is possible
for a 1-ended group $G$ to act geometrically on CAT(0) spaces
$X,Y$ such that $\bd X$, $\bd Y$ are not homeomorphic. Thus one can
not talk about `the' boundary of a CAT(0) group as in the case of
hyperbolic groups.

We remark that if a one ended hyperbolic group $G$ splits over a
two ended group $C$, then the limit points of $C$ separate the
boundary of $G$. This holds also for $CAT(0)$ groups.

We show that the converse also holds, so one can detect splittings
of $CAT(0)$ groups from their boundaries:
\begin{Thm} Let $G$ be a one ended group acting
geometrically on a CAT(0) space $X$. If  $\bd X$ has a cut pair then either $G$ splits over a 2-ended group
or $G$ is virtually a surface group.
 \end{Thm}
We remark that in case $\bd X$ is locally connected the theorem
above follows from \cite{PAP}. Indeed in this case if $\bd X$ has
a cut pair, it can be shown that a quasi-line coarsely separates
$X$.

We showed in \cite{P-S} that if $Z$ is a continuum without cut
points, then we can associate to $Z$ an $\mathbb{R}$-tree $T$
(called the JSJ tree of $Z$) encoding all pairs of points that
separate $Z$.

We show here that if $G$ is a 1-ended CAT(0) group and $Z$ is any CAT(0)
boundary then the $\mathbb{R}$-tree $T$ is simplicial and gives
the JSJ decomposition of $G$. More precisely we have:

\begin{Thm} Let $G$ be a one ended group acting
geometrically on a CAT(0) space $X$. Then the JSJ-tree of $\bd X$
is a simplicial tree $T$ and the graph of groups $T/G$ gives a
canonical JSJ decomposition of $G$ over 2-ended groups.
\end{Thm}

We note that the JSJ decomposition in the previous Theorem is
closely related to the JSJ decomposition constructed by
Scott-Swarup (\cite{S-S}).

Our strategy for obtaining these results is similar to the one of
Bowditch in the case of hyperbolic groups. For example to show
that $\bd G$ has no cut point, Bowditch associates a tree $T$ to
$\bd G$ and studies the action of $G$ on $T$ via Rips theory.

 The main difficulty in
the case of $CAT(0)$ groups is that $G$ does not act on the
boundary as a convergence group so it is not immediate, as it is
in the hyperbolic case, that the action on the trees we construct
has no global fixed point. To deal with this difficulty we show
that in the CAT(0) case the action has a convergence type
property:
\begin{Thm} Let $X$ be a CAT(0) space and $G$ a group acting
properly on $X$.  For any sequence of distinct group elements of
$G$, there exists a subsequence $(g_i)$ and points $n,p \in \bd X$
such that for any compact set $C \subset \bd X
-\bar B_T(n,\theta)$, $g_n(C) \to \bar B_T(p,  \pi
-\theta)$ (in the sense that for any open $U \supset
\bar B_T(p, \pi -\theta)$, $g_i(C) \subset U$ for all $i$
sufficiently large).
\end{Thm}

We denote  by $\bar B_T(n,\theta)$ the closed Tits ball with center $n$
and radius $\theta $. See \cite{SWE} and \cite{KAR} for partial results of this type.

The above Theorem plays a crucial role in our proofs. We think it
is also of independent interest as it is a useful tool to study
actions of CAT(0) groups on their boundaries. For example if you quotient out by  the Tits components,  $G$ acts on this quotient as a convergence group (see Corollary \ref{C:conv}).   More importantly, we obtain the
following improvement of a result of Ballmann and Buyalo \cite{BAL-BUY}:
\begin{Thm}
If the Tits diameter of $\partial X$ is bigger than $\frac {3\pi}2
$ then $G$ contains a rank 1  hyperbolic element. In particular:
if $G$ doesn't fix a point of $\bd X$ and doesn't have rank 1, and
$I$ is a minimal closed invariant set for the action of $G$ on
$\partial X$, then for any $x \in
\partial X$, $d_T(x,  I) \le \frac \pi
2$.\end{Thm}

It follows from this Theorem that if the Tits diameter of
$\partial X$ is bigger than $\frac {3\pi} 2 $ then $G$ contains a
free subgroup of rank 2 (\cite {BAL-BRI}, theorem A). We remark
that the Tits alternative is not known for CAT(0) groups. In fact
Ballmann-Buyalo (\cite{BAL-BUY}) conjecture that if $G$ acts
geometrically on $X$ and the Tits diameter of $\bd X$ is bigger
than $\pi $ then $G$ contains a rank 1 element (so $G$ contains a
free subgroup of rank 2). In the remaining case of Tits diameter
$\pi $ they conjecture that $X$ is either a symmetric space or a
Euclidean building or reducible. We note that the Tits alternative
is known in the case of symmetric spaces and Euclidean buildings,
so it would be implied for all CAT(0) groups by an affirmative answer to the above
conjectures. On the other hand if $X$ is a cube complex then the
Tits alternative holds as shown by Sageev and Wise (\cite {S-W}).
\subsection{Plan of the paper}
In section 2 we give a summary of our earlier paper \cite{P-S}. We
explain how to produce trees from cut points and cut pairs.
More precisely, we show that if $Z$ is a compact connected metric
space, then there is an $\Bbb R$-tree $T$ `encoding' all cut points
of $Z$. In particular if $Z$ has a cut point, then $T$ is non
trivial. The construction of $T$ is canonical, so the
homeomorphism group of $Z$ acts on $T$. If $Z$ has no cut points,
we show that there is an $\Bbb R$-tree $T$ `encoding' all cut
pairs of $Z$ (the `JSJ tree' of $Z$). The JSJ tree is also
canonical.

In section 3 we recall background and terminology for CAT(0)
groups and spaces. In section 4 we show that actions of CAT(0)
groups on boundaries satisfy a convergence type property ($\pi
$-convergence). Sections 5 and 6 are devoted to the proof that
CAT(0) boundaries of one ended groups have no cut points. In
section 5 we apply the machinery of section 2 to construct an
$\Bbb R$-tree on which $G$ acts. Using $\pi $-convergence we show
that the action is non-trivial. In section 6 we show that the
action is stable and we apply the Rips machine to show that $T$ is
in fact simplicial. From this we arrive at a contradiction. In
section 7 we assume that $\bd X$ has a cut pair. We apply the
construction of section 2 and we obtain the JSJ $\Bbb R$-tree for
$\bd X$. Using $\pi $-convergence again we show that the action is
non-trivial, and applying the Rips machine we show that the tree is
in fact simplicial. We deduce that $G$ splits over a 2-ended group
unless it is virtually a surface group. We show further that the
action on the tree gives the JSJ decomposition of $G$.

We would like to thank the Max Planck Institute for its
hospitality while this work was being completed in the fall of
2006.

\section{Trees and continua}

\subsection{From cut points to trees}
 Whyburn (\cite {WHY}) was the first to study the cut
point set of a locally connected metric space and to show that it
is `treelike' (a dendrite). Bowditch (\cite {BOW5}) showed how to
construct $\Bbb R$-trees for boundaries of hyperbolic groups while
Swenson (\cite {SWE}) gave a more direct construction. In \cite
{P-S} we showed how to associate to a continuum $Z$ an
$\mathbb{R}$-tree $T$ encoding the cut points of $Z$. We recall
here the results and terminology of \cite {P-S}.

In the following $Z$ will be a (metric) continuum.
\begin{Def}
If $a,b,c\in Z$ we say that $c\in (a,b)$ if $a,b$ lie in distinct
components of $Z-\{c \}$.

We call $(a,b)$ an interval and this relation an interval
relation. We define closed and half open intervals in the obvious
way i.e. $[a,b)=\{a \}\cup (a,b)$, $[a,b]=\{a,b \}\cup (a,b)$ for
$a\ne b$ and $[a,a)=\emptyset, \, [a,a]=\{a \} $.
\end{Def}
We remark that if $c\in (a,b)$ for some $a,b$ then $c$ is a cut
point. Clearly $(a,a)=\emptyset $ for all $a\in Z$.
\begin{Def} We define an equivalence relation on $Z$. Each cut point
is equivalent only to itself and if $a, b\in Z$ are not cut points
we say that $a$ is equivalent to $b$, $a\sim b$ if
$(a,b)=\emptyset $.
\end{Def}
Let's denote by $\cP$ the set of equivalence classes for this
relation. We can define an interval relation on $\cP$ as follows:
\begin{Def}
For $x,y,z\in \cP$ we say that $z\in (x,y)$ if for all $a\in
x,b\in y, c\in z$ we have that

$$[a,c)\cap (c,b]=\emptyset $$

\end{Def}

We call elements of $\cP $ which are not cut points
\textit{maximal inseparable sets}.

If $x, y,z\in \cP$ we say that $z$ is between $x,y$ if $z\in
(x,y)$. It is shown in \cite{P-S} that $\cP$ with this betweeness
relation is a pretree.

The $\mathbb{R}$ tree $T$ is obtained from $\cP$ by ``connecting
the dots" according to the pretree relation on $\cP$. We proceed
now to give a rigorous definition of $T$.

We have the following results about intervals in pretrees from
\cite{BOW5}:

\begin{Lem} \label{L:subset}
If $x,y,z \in \cP$, with $y\in [x,z]$ then $[x,y]\subset [x,z]$.
\end{Lem}

\begin{Lem} \label{L:order}
Let $[x,y]$ be an interval of $\cP$. The interval structure
induces two linear orderings on $[x,y]$, one being the opposite of
the other, with the property that if $<$ is one of the orderings,
then for any $z,w\in [x,y]$ with $z<w$, $(z,w)=\{ u\in [x,y]:z<u<w
\}$. In other words the interval structure defined by one of the
orderings is the same as our original interval structure.
\end{Lem}
\begin{Def} If $x,y$ are distinct points of $\cP$ we say that $x,y$
are \textit{adjacent} if $(x,y)=\emptyset $. We say $x \in \cP$ is
{\em terminal} if there is no pair $y,z \in \cP$ with $x \in
(y,z)$.
\end{Def}
We recall  \cite[Lemma 4]{SWE}.
\begin{Lem} \label{L:adjacent}
If $x,y \in \cP$, are adjacent then exactly one of them is a cut
point and the other is a nonsingleton equivalence class whose
closure contains this cut point.
\end{Lem}
We have the Theorem  \cite[Theorem 6]{SWE}:
\begin{Thm}\label{T:nested}
A nested union of intervals of $\cP$ is an interval of $\cP$.
\end{Thm}
\begin{Cor}\label{C:supremum} Any interval of $\cP$ has the supremum
property with respect to either of the linear orderings derived
from the interval structure.
\end{Cor}
\begin{proof}
Let $[x,y]$ be an interval of $\cP$ with the linear order $\leq $.
Let $A\subset [x,y]$. The set $\{[x,a]:a \in A\}$ is a set of
nested intervals so their union is an interval $[x,s]$ or $[x,s)$
and $s=sup\,A$.
\end{proof}
\begin{Def} A {\em big arc} is the homeomorphic image of a compact connected
nonsingleton linearly ordered topological space.  A separable big
arc is called an {\em arc}. A {\em big tree} is a  uniquely
big-arcwise connected topological space.  If all the big arcs of a
big tree are arcs, then the big tree is called a {\em real tree}.
A metrizable real tree is called an $\R$-tree. An example of a
real tree which is not an $\mathbb{R}$-tree is the long line (see
\cite{HOC-YOU} sec.2.5, p.56).
\end{Def}

\begin{Def} A pretree $\cR$ is {\em complete} if  every closed interval  is
complete as a linearly ordered topological space (this is slightly
weaker than the definition given in \cite{BOW5}). Recall that a
linearly ordered topological space is complete if every bounded
set has a supremum.

Let $\cR$ be a pretree, an interval $I\subset \cR$ is called {\em
preseparable} if there is a countable set $Q \subset I$ such that
for every nonsingleton closed interval $J \subset I$, $J \cap Q
\neq \emptyset$. A pretree is {\em preseparable} if every interval
in it is preseparable.
\end{Def}
We recall now the construction in \cite{P-S}.  Let $\cR$ be a
complete pretree. Set
$$T=\cR \cup \underset {{x,y\,\, adjacent}} { \bigsqcup } I_{x,y}$$
where $I_{x,y}$ is a copy of the real open interval $(0,1)$ glued
in between $x$ and $y$. We extend the interval relation of $\cR$
to $T$ in the obvious way (as in \cite{SWE}), so that in $T$,
$(x,y) = I_{x,y}$.    It is clear that $T$ is a complete pretree
with no adjacent elements.   When $\cR=\cP$, we call the $T$ so
constructed the cut point tree of $Z$.

\begin{Def} For $A$ finite subset of $T$ and $s\in T$ we define
$$U(s,A)=\{ t\in T:[s,t]\cap A=\emptyset \}$$
\end{Def}

The following is what the proof of \cite[Theorem 7]{SWE}   proves
in this setting.
\begin{Thm} Let $\cR$ be a complete pretree.  The pretree $T$, defined
above, with the topology defined by the basis  $\{U(s,A)\}$ is a
regular big tree.  If in addition $\cR$ is preseparable, then $T$
is a real tree.
\end{Thm}

We recall  \cite[Theorem 13]{P-S}
\begin{Thm} \label{T:countable}
The pretree $\cP$ is preseparable, so the cut point tree $T$ of
$Z$ is a real tree.
\end{Thm}

The real tree $T$ is not always metrizable.  Take for example $Z$
to be the cone on a Cantor set $C$ (the so called Cantor fan).
Then $Z$ has only one cut point, the cone point $p$, and $\cP $
has uncountable many other elements $q_c$, one for each point $c
\in C$. As a pretree, $T$ consists of uncountable many arcs $\{
[p,q_c]:\, c \in C\}$ radiating from a single central point $p$.
However, in the topology defined from the basis $\{U(s,A)\}$,
every open set containing $p$ contains the arc $[p,q_c]$ for all
but finitely many $c \in C$. There can be no metric, $d$, giving
this topology since $d(p,q_c)$ could only be non-zero for
countably many $c \in C$. In \cite {P-S} we showed that it is
possible to equip $T$ with a metric that preserves the pretree
structure of $T$. This metric is `canonical' in the sense that any
homeomorphism of $Z$ induces a homeomorphism of $T$. We recall
briefly how this is done. The idea is to metrize $T$ in two steps.
In the first step one metrizes the subtree obtained by the span of
cut points of $\cP$. This can be written as a countable union of
intervals and it is easy to equip with a metric.

$T$ is obtained from this tree by gluing intervals to some points
of $T$. In this step, one might glue uncountably many intervals, but
the situation is similar to the Cantor fan described above. The
new intervals are metrized in the obvious way, e.g. one can give
all of them length one.

So we have the following:
\begin{Thm}\label{T:Rtree}
There is a path metric $d$ on $T$, which preserves the pretree
structure of $T$, such that $(T,d)$ is a metric $\mathbb{R}
$-tree. The topology so defined on $T$ is  canonical (and may be
different from the topology with basis $\{U(s,A)\}$). Any
homeomorphism $\phi $ of $Z$ induces a homeomorphism $\hat \phi $
of $T$ equipped with this metric.
\end{Thm}

Since $T$ has the supremum property, all ends of $T$ correspond to elements of $\cP$, and these elements will be terminal in $\cP$ and in $T$.  We remove these terminal points to obtain a new
tree which we still call $T$. Clearly the previous Theorem holds
for this new tree too.

\subsection{JSJ trees for continua }
In \cite{P-S} we showed how to associate to cut-pairs of a
continuum an $\Bbb R$-tree which we call the JSJ-tree of the
continuum. The construction is similar as in the case of cut
points but for JSJ-trees a new type of vertices has to be
introduced which corresponds to the hanging orbifold vertex groups
in the group theoretic setting. We recall briefly now the results
of \cite{P-S}.

 Let $Z$ be a metric continuum without cut points.  The set
$\{a,b\}$ is a cut pair of $Z$ separating $p\in Z$ from $q \in Z$
if there are continua $P\ni p$ and $Q\ni q$ such that $P \cap Q
=\{a,b\}$ and $P \cup Q = Z$.
\begin{Def} Let $Z$ be a continuum without cut points.
A finite set $S$ with $|S|>1$ is called a {\em cyclic} subset if  $S$ is a
cut pair, or if there is an ordering $S=\{x_1,\dots x_n\}$ and
continua $M_1,\dots M_n$ with the following properties:
\begin{itemize}
\item  $ M_n \cap M_1 = \{x_1\} $,  and for $i>1$, $\{x_i\} =
M_{i-1} \cap M_{i}$ \item $M_i \cap M_j = \emptyset$ for $i-j \neq
\pm 1 \mod n$ \item $\bigcup M_i = Z$
\end{itemize}
The collection $M_1, \dots M_n$ is called the (a) {\em cyclic
decomposition} of $Z$ by $\{x_1,\dots x_n\}$. (For $n>2$, this
decomposition is unique.) If $S$ is an infinite subset of $Z$ and
every finite subset $A \subset S$ with $|A|>1$ is cyclic, then we
say $S$ is cyclic.  Every cyclic set is contained in a maximal
cyclic set.   We call a maximal cyclic subset of $Z$ with more than 2 points a {\em necklace}.
\end{Def}
Let $S$ be a cyclic subset of $Z$. There exists a continuous
function $f: Z \to S^1$ with the following properties:
\begin{enumerate}
\item The function $f$ is one to one on $S$, and $f^{-1}(f(S))=S$
\item For $x,y\in Z$  and  $a,b\in S$ :
\begin{enumerate}
\item If  $\{f(a), f(b)\}$ separates $f(x)$ from $f(y)$  then
$\{a,b\}$ separates $x$ from $y$. \item  If $x \in S$ and
$\{a,b\}$ separates $x$ from $y$, then $\{f(a), f(b)\}$ separates
$f(x)$ from $f(y)$
\end{enumerate}
\end{enumerate}
Furthermore if $|S|>2$ (if $S$ is a necklace for example), this function can be chosen in a canonical
fashion up to isotopy and orientation.  This function will be
called the circle function for $S$ and denoted $f_S$.
\begin{Def} Let $Z$ be a metric space without cut points.
A non-empty non-degenerate set $A \subset Z$ is called inseparable
if no pair of points in $A$ can be separated by a cut pair.
  Every inseparable set is contained in a maximal inseparable set.

\end{Def}
 We define $\cP$ to be the collection of all necklaces of $Z$,
  all maximal inseparable subsets of $Z$, and all inseparable cut pairs of
$Z$.

 If two elements of $\cP$ intersect then the intersection is either a single
point, or an inseparable cut pair.
  There is a natural pre-tree structure on $\cP$ (given by separation).
   The pre-tree $\cP$ has the property that every monotone sequence in
    a closed interval of $\cP$ converges in that interval.

  Distinct points of $\cP$ are adjacent if there there are no points of
$\cP$ between them.
  \begin{Lem}\cite[Lemma 28]{P-S} If $A, B\in \cP$, with $A \subset B$ (thus $A$ is an inseparable pair and
$B$ is maximal cyclic or maximal inseparable) then $A$ and $B$ are
adjacent. The only other way two points of $P$ can be adjacent is
if one of them, say $A$, is a necklace and the other, $B$, is
maximal inseparable  and $|\bd A \cap B| \ge 2$.
\end{Lem}
  The Warsaw circle gives an example of this.

Recall that a point of $\cP$ is called terminal if it is not in any open interval of $\cP$.
Terminal points of $\cP$ do arise in CAT(0) boundaries even when the action on $\cP$ is nontrivial (see \cite{C-K} for example).
We first remove all terminal points  from  $\cP$, and then glue
open intervals between the remaining adjacent points of  $\cP$ to form a
topological $\R$-tree $T$.
  We call $T$ the
  \textbf{JSJ-Tree} of the continuum $Z$.  Notice that in the case where $Z$ is a circle,
   $\cP$ has only one element, a necklace, and so $T$ is empty.
    This and similar issues will be dealt with by reducing to the case where $G$ doesn't fix an element of $\cP$, in which case $T$ will be non-trivial.

\section{CAT(0) Groups and boundaries}
We give a brief introduction to CAT(0) spaces (see \cite{BRI-HAE},
\cite{BAL-BRI} for details). Let $Y$ be a  metric space, and $I$
be an interval of $\R$. A path $\gamma:I\to X$  is called a {\em
geodesic} if for any $[a,b] \subset I$, $$\ell_a^b (\gamma) =
b-a=d(\gamma(a),\gamma(b))$$ where $\ell_a^b(\gamma)$ is the
length of $\gamma$ from $a$ to $b$, defined as
$$\ell_a^b(\gamma)=\sup\left \{\sum\limits_{i=1}^n
d\left(\gamma(x_{i-1}),\gamma(x_i)\right) : \, \{x_0, \dots x_n\}
\text{ is a partition of }[a,b] \right\}$$ If $I$ is a ray or $\R$
then we refer to $\gamma$ (or more precisely its image) as a
geodesic ray or geodesic line respectively.  If $I$ is a closed
interval, we refer to $\gamma$ as a geodesic segment.

A metric space is called {\em proper}  if every closed metric ball
is compact. In a proper metric space, closed and bounded implies
compact.

A CAT(0) space $X$ is a proper geodesic metric space with the
property that every geodesic triangle (the union of three geodesic
segments meeting in the obvious way) in $X$ is a least as thin as
the comparison Euclidean triangle, the triangle in the Euclidean
plane with the same edge lengths.   Every CAT(0) space is
contractible.

The (visual) boundary $\bd X$ is the set of equivalence classes of
unit speed geodesic rays, where $R,S:[0,\iy) \to X$ are equivalent
if $d(R(t), S(t))$ is bounded.  Given a ray $R$ and a point $x \in
X$ there is a ray $S$ emanating from $x$ with $R \sim S$.  Fixing
a base point $\mathbf 0 \in X$ we define a Topology on $\bar X= X
\cup \bd X$ by taking the basic open sets of $ x \in X$ to be the
open metric balls about $x$. For $y \in \bd X$, and $R$ a ray from
$\mathbf 0$ representing $y$, we construct  basic open sets
$U(R,n,\epsilon)$ where $n,\epsilon>0$. We say $z \in
U(R,n,\epsilon)$ if the unit speed geodesic, $S:[0,d(\mathbf 0,z)]
\to \bar X$, from $\mathbf 0$ to $z$ satisfies $d(R(n),S(n))
<\epsilon$.  These sets form a basis for a topology on $\bar X$
under which $\bar X$ and $\bd X$ are compact metrizable.
Isometries of $X$ act by homeomorphisms on $\bd X$.

For three points $a,b,p \in X$ $a\neq p\neq b$ the {\em comparison
angle} $\overline{\angle}_p(a,b)$ is the Euclidean angle at the
point corresponding to $p$ in the comparison triangle in $\E^2$.
By CAT(0), for any $c\in (p,a)$ and $d \in (p,b)$,
$\overline{\angle}_p(a,b)\ge \overline{\angle}_p(c,d)$.  This
monotonicity implies that we can take limits.  Thus the {\em
angle} is defined as  $\angle_p(a,b)  = \lim
\overline{\angle}_p(a_n,b_n)$  where $(a_n) \subset (p, a)$ with
$a_n \to p$ and $(b_n) \subset (p, b)$ with $b_n \to p$. By
monotonicity this limit exists, and this definition works equally
well if one or both of $a,b$ are in $\bd X$. Similarly we can
define the {\em comparison angle}, $\overline{\angle}_p(a,b)$ even
if one or both of $a,b$ are in $\bd X$ (the limit will exist by
monotonicity).  By \cite[II 1.7(4)]{BRI-HAE},
$\overline{\angle}_p(a,b) \ge \angle_p(a,b)$.

On the other hand for $a,b \in \bd X$, we define $\angle(a,b) =
\sup\limits_{p \in X} \angle_p(a,b)$.
 It follows from \cite[II 9.8(1)]{BRI-HAE} that $\angle(a,b) =
\overline{\angle}_p(a,b)$ for any $p \in X$.
 Notice that isometries of $X$ preserve the angle between points of $\bd X$.
 The angle defines a path metric, $d_T$ on the set $\bd X$, called the Tits
metric, whose topology is finer than the given topology of $\bd
X$.  Also $\angle(a,b)$ and $ d_T(a,b)$ are equal whenever either
of them are less than $\pi$ (\cite[II 9.21(2)]{BRI-HAE}.

  The set $\bd X$ with the Tits metric is called the Tits  boundary of $X$,
denoted  $TX$.  Isometries of $X$ extend to isometries of $TX$.

  The identity function $TX \to \bd X$ is continuous (\cite[II 9.7(1)]{BRI-HAE}), but the identity
function
 $\bd X \to TX$ is only lower semi-continuous(\cite[II 9.5(2)]{BRI-HAE}).  That is for any sequences
$(y_n), (x_n)  \subset \bd X$ with $x_n \to x$ and $y_n \to y$ in
$\bd X$, then $$\varliminf d_T(x_n,y_n) \ge d_T(x,y)$$

Recall that the action of a group $H$ on a space $Y$ is called:
\begin{itemize}
\item {\em Proper} if
 for every compact $K \subset Y$, the set $\{h \in H : h(K)
\cap K \neq \emptyset \}$ is finite.
\item {\em Cocompact} if the quotient space of $Y$ by the
action of $H$,  $Y /H$ is compact.
\item {\em Geometric} if $H$ acts properly,
cocompactly by isometries on $Y$.
\end{itemize}
It follows that $H$ is
quasi-isometric to $Y$, that is they have the same coarse
geometry.

A {\em hyperbolic} isometry $g$ of a CAT(0) space $X$ is an
isometry which acts  by translation on a geodesic line $L\subset
X$.  The distance that $g$ translates along $L$ is called the {\em
translation length} of $g$ denoted $|g|$. Any such  line $L$ is
called an  {\em axis} of $g$ and the endpoints of $L$ (in $\bd X$)
are denoted by $g^+$ and $g^-$, where $g^+$ is the endpoint in the
direction of translation. If a group acts geometrically on a
CAT(0) space, then the non-torsion elements are exactly the
hyperbolic elements, and every group element $g$ has a {\em
minimum} set  ($\Min g$) which is the set of points moved a
minimum distance by $g$.  When $g$ is hyperbolic, the $\Min g$ is
just the union of the axes of $g$.  When $g$ is torsion (called
{\em elliptic}) then  $\Min g$ is just the set of points of $X$
fixed by $g$. In either case, $\Min g$ is convex and the
centralizer $Z_g$ acts geometrically on $\Min g$.

If a group $G$ acts geometrically on a non-compact CAT(0) space
$X$, then $G$ has a hyperbolic  element, and $\bd X$ is
finite dimensional (\cite[Theorem 12]{SWE}).
\begin{Def} If $X$ is a CAT(0) space and $A\subset X$ we define
$\Lambda A$ to be the set of limit points of $A$ in $\bd X$. If $H$  is a group acting properly by isometries on the CAT(0) space $X$  we define $\Lambda
H=\Lambda Hx$ for some $x\in X$. Observe that $\Lambda H$ does not
depend on the choice of $x$.
\end{Def}
{\bf For the remainder of the paper, $G$ will be a group
acting geometrically on the non-compact CAT(0) space $X$}

\begin{Lem}  If $h\in G$ is hyperbolic and the centralizer $Z_{h^n}$ is virtually cyclic
for all $n$, then for some $n$,  $stab(h^\pm)= Z_{h^n}$
\end{Lem}
\begin{proof}
By \cite[Zipper Lemma]{SWE}  hyperbolic elements with the same endpoints have a
common power. Thus we may assume that $h$ has the minimal
translation length of any hyperbolic element with endpoints
$h^\pm$.
 Let $H= stab(h^\pm)$. Let $Z$ be the union of all lines from $h^+$ to
$h^-$.  By \cite{BRI-HAE} $W = R\times Y$
 where $R$ is an axis of $h$ and $Y$ is a convex subset of $X$, and $
stab(h^\pm)$ acts on $W$ preserving the product structure.
 Since $h$ has minimal translation length,  every element of $H$
translates the $R$ factor by a multiple of $|h|$.

By \cite[Theorem 8]{SWE} $stab(h^\pm) = \bigcup\limits_n Z_{h^n}$.   Since
$Z_{h^n} \cup Z_{h^m} < Z_{h^{nm}}$, it suffices to uniformly bound the
number of cosets of $\langle h \rangle$ in $Z_{h^n}$.  We do this by
finding a section of the projection function $Z_{h^n} \to Z_{h^n}/ \langle h
\rangle$ whose image generates a finite group. Since there are only
finitely many conjugacy classes of finite subgroups of $G$, there is a uniform bound on the size of a finite subgroup of $G$, and this number also bounds the number of cosets of $\langle h\rangle $ in  $Z_{h^n}$ for all $n$.

Since $Z_{h^n}$ is virtually cyclic, the number of cosets  of $\langle h \rangle$ in $Z_{h^n}$  is finite (say $p$).  Choose coset representatives  $g_1, \dots g_p \in  Z_{h^n}$.   Multiplying by a multiple of $h$,
we may assume that each  $g_i$ fixes the $R$ factor of $W$.  Thus
$K=\langle g_1,g_2,\dots g_p\rangle $ is a subgroup of the
virtually cyclic $Z_{h^n}$ which doesn't translate along the $R$
direction. It follows that no power of $h$ lies in $K$.  If two subgroups of a virtually $\Z$ group intersect trivially then one of them is finite, and so  $K$
is finite.
\end{proof}
\begin{Lem}\label{L:limsets} Let $K$ be a group acting geometrically on a CAT(0) space $U$. If $H<K$ with $\Lambda H$  finite then $H$ is virtually cyclic.
\end{Lem}
\begin{proof}
We assume that $H$ is infinite.  Since $\Lambda H$ is finite, passing to a  finite index  subgroup, we may assume that $H$ fixes $\Lambda H$ pointwise.

We first show that $H$ contains a hyperbolic element. Suppose not, then from \cite[Theorem 17]{SWE} there is a convex set $Y \subset U$ with $\dim \bd Y <\dim \bd U$
such that $Y$ is invariant under the action of $H$ and the action of $H$ on $Y$ can be extended to a geometric action on $Y$.  Since $\dim \bd U <\iy$ (\cite{SWE}),
 we get a contradiction by induction on boundary dimension (In the case where the boundary is zero dimensional, any group acting geometrically is virtually free.)

Thus there exists $h \in H$ hyperbolic. By \cite{SWE} $stab(h^\pm)
= \bigcup\limits_n Z_{h^n}$, and so $H < \bigcup\limits_n
Z_{h^n}$.   The union $W$ of axis of $h$ from $h^-$ to $h^+$
decomposes as $W = R \times Y$ where $R$ is an axis for $h$ and
$Y$ is convex in $U$. Since $H$ fixes  $\{h^\pm\}$,  then $H$ acts
on $W$ preserving this product structure.

Suppose that $H/\langle h\rangle$ is infinite.  Let $\epsilon$ be the translation length of $h$.  Choose one element
$g_i$ from each coset with the translation distance of $g_i$ on
the $R$ factor less than $\epsilon$.  Choose $m$ larger than $|\{g \in K \,|\, g(B(u,\epsilon)) \cap B(u, \epsilon) \neq \emptyset \}|$
 for any $u \in U$.   Choose $n$ such that
$g_1, \dots g_m \in Z_{h^n}$.   Notice that  $\Min h^n= R\times Z$ where $Z$ is a convex subset of $Y$.
  Now consider the images $\hat g_1, \dots \hat g_m$ of $g_1,\dots g_m$ in the quotient  $Z_{h^n}/\langle h^n\rangle=J$ which acts geometrically on $Z$.
 Notice that  $\hat g_1, \dots \hat g_m$ are  distinct elements of $J$.
  By choice of $m$  the subgroup $F = \langle \hat g_1, \dots \hat g_m \rangle$ doesn't fix a point of $Z$.  It follows that $F$ is infinite.
   Notice that $\Lambda F \subset \bd Z$ is finite (in fact $\Lambda F \subset \Lambda H$).  It follows by the above that $F$ contains a hyperbolic element $\hat g$.
    Let $g$ be a preimage of $\hat g$ in $H$.  Clearly $g$ is hyperbolic and $\{g^\pm\} \cap \{h^\pm\} = \emptyset$.  However $g \in Z_{h^j}$ for some $j$, so
$\langle g, h^j\rangle \cong \Z^2$ and has limit set a circle.  This contradicts the fact that $\Lambda H$ is finite.

\end{proof}

\section {Convergence actions}
\subsection{$ \pi $-Convergence Action}
  Notice that $\overline{\angle}_p(a,b)$ is a function of
$d(a,b)$, $d(p,a)$ and $d(b,p)$.  It follows that the function
$\overline{\angle} $ is continuous on the subset of $X^3$ ,
$\{(a,b,p)\in X^3 : a\neq p,\, b \neq p\}$.

\begin{Lem}\label{L:atmostpi} Let $p,q \in X$ $p \neq q$ and $a \in \bd X$.
$\overline{\angle}_p(q,a)+ \overline{\angle}_q(p,a) \le \pi$.
\end{Lem}
\begin{proof} Suppose not, then by the definition of these angles as monotone
limits, there exists $p' \in [p,a)$ and $q' \in [q,a)$ such that
$\overline{\angle}_p(q,p')+ \overline{\angle}_q(p,q') > \pi$.  We
may also assume that $d(p,p') = d(q,q')$.  By convexity of the
metric applied to the rays $[p,a)$ and $[q,a)$ we have that
$d(p',q') \le d(p,q)$.

By \cite[II 1.1]{BRI-HAE}, for any quadrilateral (in our case
$pp'q'q$) in $X$ there is a comparison quadrilateral ($\bar p \bar
p' \bar q' \bar q$) in $\E^2$ having the same edge lengths
($d(p,p') = d(\bar p,\bar p')$, $d(p',q') = d(\bar p',\bar q')$
,$d(q,q') = d(\bar q,\bar q')$ and $d(q,p)= d(\bar q,\bar p)$)
with diagonals no shorter than in the original quadrilateral
($d(p,q') \le d(\bar p,\bar q')$ and $d(q,p') \le d(\bar q,\bar
p')$ ). It follows that $\overline{\angle}_q(p,q') \le
\angle_{\bar q}(\bar p, \bar q') $ and that
$\overline{\angle}_{p}(q,p') \le \angle_{\bar p}(\bar q, \bar p')
$.   Thus $ \angle_{\bar q}(\bar p, \bar q') + \angle_{\bar
p}(\bar q, \bar p') > \pi$.  Using high school geometry we see
that $d(\bar p', \bar q') > d(\bar p, \bar q)$  which contradicts
the fact that $d( p',  q') \le d( p,  q)$.
\end{proof}
\begin{Lem}\label{L:compproj} Let $X$ be a CAT(0) space,
$\theta \in (0, \pi)$, $q\in \bd X$, and $K$ a compact subset of
$\bd X- \bar B_T(q, \theta)$.  Then for any $x \in X$, there is a
point $y \in [x,q)$ and $\epsilon >0$ with $\angle_y(q,c)
>\theta +\epsilon$ for all $c \in K$.
\end{Lem}
\begin{proof}  We first find $y \in [x,q) $ such that $\angle_y(q,c)
>\theta $ for all $c \in K$.
Suppose there is not such $y$, then by \cite[II 9.8 (2)]{BRI-HAE},
there is a sequence of points $(c_i)\subset K$ and  a monotone
sequence $(y_i) \subset [x,q)$ such that $y_i \to q$ with
$\angle_{y_i}(q,c_i)\le \theta $.  By the monotonicity of $(y_i)$
and \cite[II 9.8 (2)]{BRI-HAE}
$$\angle_{y_i}(q,c_n)\le \theta$$ for all $n \ge i$.  By
compactness of $K$, passing to a subsequence we may assume that
$c_i \to c \in K$.  For any fixed $i \in \N$, the sequence of rays
$\left([y_i, c_n)\right)$ converges to the ray $[y_i,c)$.  Thus by
\cite[II 9.2(1)]{BRI-HAE} $\angle_{y_i}(q,c_n) \to
\angle_{y_i}(q,c) \le \theta$ and this is true for each $i \in\N$.
It follows from \cite[II 9.8 (2)]{BRI-HAE} that $\angle(q,c) \le
\theta$ which is a contradiction since $c\not \in
\overline{B}_T(q, \theta)$.  Thus we have  $y \in [x,q) $ such
that $\angle_y(q,c)
>\theta $ for all $c \in K$.  Since the function $\angle_y$ is
continuous on  its domain, by compactness of $K$, there is a $d
\in K$ with $\angle_y(q,d) \le \angle_y(q,c)$ for all $c \in K$.
We let $\epsilon <\theta -\angle_y(q,d)$ and the result follows.
\end{proof}
\begin{Lem}\label{L:pi}
Let $X$ be a CAT(0) space and $G$ a group acting properly on $X$.
Let $x \in X$, $\theta \in [0, \pi]$ and $(g_i)\subset G$ with the
property that $g_i(x) \to p \in \bd X$ and $g_i^{-1}(x) \to n \in
\bd X$. For any compact set $K \subset \bd X
-\overline{B}_T(n,\theta)$, $g_n(K) \to \overline{B}_T(p,  \pi
-\theta)$ (in the sense that for any open $U \supset
\overline{B}_T(p, \pi -\theta)$, $g_n(K) \subset U$ for all $n$
sufficiently large).
\end{Lem}
The case where $\theta= \frac \pi 2$ was done independently by
Karlsson \cite{KAR}.
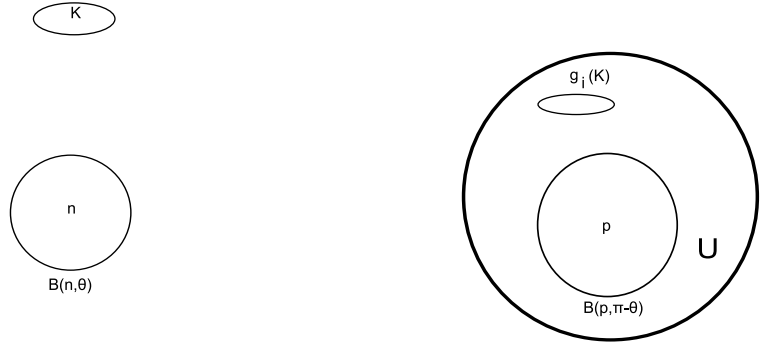
\begin{figure}[h]
\includegraphics[width=4.2in ]{Pic1}
\end{figure}
\begin{proof}
It suffices to consider the case $\theta \in (0, \pi)$.

 Let $K \subset \bd X - \overline{B}_T(n, \theta)$ be compact,
and $U$ an open set in $\bd X$ containing the closed ball
$\overline{B}(p, \pi -\theta)$.

Fix $x \in X$. By Lemma \ref{L:compproj}, there is a $y \in [x,
n)$ and $\epsilon >0$ with $\angle_y(n,c) >\theta +\epsilon$ for
all $c \in K$.

Notice that the sequence segments $\left( [x,g_i^{-1}(x)]\right)$
converges to the ray $[x,n)$.  For each $i\in \N$ let $y_i
$ be the projection  of $y$ onto $[x,g_i^{-1}(x)]$.  Passing to a
subsequence, we may assume that $d(y_i,g_i^{-1}(x))\ge 1$.  Let
$z_i \in [y_i,g_i^{-1}(x)]$ such that $d(y_i,z_i) =1$.  It follows
that $z_i \to z \in [y,n)$ where $d(y,z)=1$.

\begin{figure}[htp]
\includegraphics[width=4.2in ]{Pic2}
\end{figure}

Now comparison angles are continuous in all three variables for
points on $X$ and lower semi-continuous for points on $\bd X$.
That is if $D= \{ (u,v,w) \in \overline{X}^2 \times X |\, u \neq w
\neq v\}$ then the comparison function $\overline{\angle}: D \to
[0,\pi]$ is continuous on the third variable, it is continuous on
$D'= \{ (u,v,w) \in X ^3 |\, u \neq w \neq v\}$ and it is lower
semicontinuous on the first 2 variables. This is not true for the
function $\angle: D \to [0,\pi]$, as it is only upper
semicontinuous. However  it is always that case that
$\overline{\angle}_w(u,v) \ge \angle_w(u,v)$.  Thus for any $c \in
K$,
$$\overline{\angle}_{y_i}(z_i,c) \to
\overline{\angle}_{y}(z,c)\ge \angle_{y}(z,c) = \angle_{y}(n,c)
>\theta +\epsilon$$  It follows that  for $i\gg 0$,
$\overline{\angle}_{y_i}(z_i,c) > \theta + \epsilon$ for all $c
\in K$.

Now consider the segments $g_i([x,g_i^{-1}(x)]) =[g_i(x),x]$  and
note that $ [g_i(x),x] \to (p,x]$ the ray.  Notice that for any $c
\in K$ and $i \gg 0$,
$\overline{\angle}_{g_i(y_i)}(g_i(z_i),g_i(c))=\overline{\angle}_{y_i}(z_i,c)
>\theta + \epsilon$.

Let $u \in (p,x]$, for each $i$, let $u_i$ be the projection of
$u$ to $[g_i(x),x]$.  For $i \gg 0$, $u_i \in [g_i(z_i), x]$ and
$d(u_i,g_i(z_i))>1$.  For $i\gg 0$ choose $v_i \in [g_i(z_i),u_i]$
with $d(v_i,u_i)=1$.  By  Lemma \ref{L:atmostpi} and the
monotonicity of the comparison angle,
$$\overline{\angle}_{g_i(y_i)}(g_i(z_i),g_i(c))
+\overline{\angle}_{u_i}(v_i,g_i(c)) \le \pi$$ for all $c \in K$
and $i \gg 0$.   Thus for $i \gg 0$
$$\overline{\angle}_{u_i}(v_i,g_i(c)) < \pi-\theta-\epsilon$$ for
all $c \in K$.

Notice that $v_i \to v \in (p, u]$ where $d(v,u)=1$.  Using  the
continuity of $\overline{\angle}$, we can show that for $i\gg 0$,
$ \pi -\theta> \overline{\angle}_{u}(v,g_i(c)) \ge
\angle_u(p,g_i(c)) $.  So for any $u \in (p,x]$ there is $I \in
\N$ such that for all $i>I$ and for all $c \in K$,
$\angle_u(p,g_i(c))<  \pi -\theta$.
\subsubsection*{Claim: $g_i(K) \subset U$ for all $i \gg 0$}
Suppose not, then passing to a subsequence,  for each $i$, there
exists $c_i\in K$, such that $g_i(c_i) \not \in U$.  Passing to a
subsequence we may assume that $g_i(c_i) \to \hat c \in \bd X -
B_T(p,  \pi -\theta)$.  Since $\angle(\hat c, p) >  \pi-\theta$,
there exists $w \in (p,x]$ such that $\angle_w(\hat c,p) > \pi
-\theta$. However this contradicts the fact that for $i\gg 0$,
$\angle_w(p,g_i(c_i))< \pi -\theta$.  This proves the claim and
 the Lemma.
\end{proof}
\begin{Thm}\label{T:pi} Let $X$ be a CAT(0) space and $G$ a group acting
properly  by isometries  on $X$.  For any sequence of distinct group elements of
$G$, there exists a subsequence $(g_i)$ and points $n,p \in \bd X$
such that for any $\theta \in [0,\pi]$ and for any compact set $K
\subset \bd X -\overline{B}_T(n,\theta)$, $g_n(K) \to
\overline{B}_T(p, \pi -\theta)$ (in the sense that for any open $U
\supset \overline{B}_T(p, \pi -\theta)$, $g_i(K) \subset U$ for
all $i$ sufficiently large).
\end{Thm}
\begin{proof}
 For any sequence of distinct elements of $G$, using the
compactness of $\overline{X}= X \cup \bd X$, there exists a
subsequence $(g_i)$ and points $n,p \in \bd X$ such that $g_i(x)
\to p$, $g_i^{-1}(x) \to n$ for some $x \in X$, and so Lemma \ref{L:pi} applies.
\end{proof}
We will refer to the results of Lemma \ref{L:pi} and Theorem \ref{T:pi} collectively as
{\em $\pi$-convergence}.
\begin{Cor} \label{C:conv} Let $G$ act properly by isometries on $X$, a
CAT(0) space.  If $\tau: \bd X \to Z$ is a $G$-quotient map, with the property that for any
$ a, b \in \bd X$ with $d_T(a,b) \le \frac \pi 2$, $\tau(a) =
\tau(b)$, then the induced action of $G$ on $Z$ is a convergence
action.
\end{Cor}
\begin{proof} Let $(g_i)$ be a sequence of distinct elements of $G$.
By $\pi$-convergence, passing to a subsequence there exists $n,p
\in \bd X$ such that for any compact $K \subset \bd X- B_T(n,
\frac \pi 2)$ and any open $U \supset B_T(p, \frac \pi 2)$,
$g_i(K) \subset U$ for all $i \gg 0$.  Now let $\hat p = \tau(p)$
and $\hat n = \tau (n)$.  For any compact $\hat K \subset Z$ with
$\hat n \not \in \hat K$, and any open $\hat U \ni \hat p$,
$\tau^{-1}(\hat U) = U$ is an open set containing $B_T(p,\frac \pi
2)$ and $\tau^{-1}(\hat K) = K$ is a compact set missing $B_T(n ,
\frac \pi 2)$.  Thus for $i\gg 0$, $g_i(K) \subset U$ and so
$g_i(\hat K) \subset \hat U$ as required.

\end{proof}

In the later sections we will be using the following result of
Ballmann-Buyalo (\cite[Proposition 1.10]{BAL-BUY})
\begin{Thm} \label{T:Bal}If the Tits diameter of $\partial X$ is bigger
than $2\pi $ then $G$ contains a rank 1 hyperbolic element. In
particular: If $G$ doesn't have rank 1, and $I$ is a minimal
closed invariant set for the action of $G$ on $\partial X$, then
for any $x \in
\partial X$ and $ m \in I$, $d_T(x,m)\leq \pi $.
\end{Thm}

\begin{Def}We say a hyperbolic element $g$ has {\em rank 1}, if
no axis of $g$ bounds a half flat (isometric to $[0,\iy) \times
\R$).  In particular if $d_T(g^+, g^{-}) > \pi$ then the
hyperbolic element $g$ is rank 1.  If $G$ contains a rank 1
hyperbolic element, then we say $G$ has {\em rank 1}. The endpoints of a
rank 1 hyperbolic element have infinite Tits distance from any
other point in $\bd X$  \cite{BAL}.  (In the case where $d_T(g^+, g^{-}) > \pi$, this follows trivially from $\pi$-convergence)
\end{Def}

Now if $G$ has rank  1, then for any open $U$, $V$ in $\bd X$,
there is a rank 1 hyperbolic element $g$ with one endpoint in
$U$ and the other in $V$ satisfying $g(\bd X-U) \subset V$
and $g^{-1}(\bd X -V) \subset U$ \cite[Theorem A]{BAL-BRI}  (Notice that $U$ and $V$ need not be disjoint.)
If $G$ has rank
1, then the only nonempty closed invariant subset of $\bd X $ is $\bd
X$.
\begin{Def}
Let $Z$ be a metric space and $A \subset Z$.  The radius of $A$, $$radius (A) = \inf\{r :\, A \subset \bar B(z,r) \text{ for some } z \in Z\}$$
We say that $c \in Z$ is a {\em centroid} of $A \subset Z$ if $A \subset \bar B(c, radius A)$.
When using the Tits metric on $\bd X$ we signal this to the reader using the subscript $_T$.
\end{Def}

Using $\pi $-convergence, we obtain the following strengthening of
the Theorem above.
\begin{Thm} \label{T:Ball}

If the Tits diameter of $\partial X$ is bigger than $\frac {3\pi}2
$ then $G$ contains a rank 1  hyperbolic element. In particular:
if $G$ doesn't fix a point of $\bd X$ and doesn't have rank 1, and
$I$ is a minimal closed invariant set for the action of $G$ on
$\partial X$, then for any $x \in
\partial X$, $d_T(x,  I) \le \frac \pi
2$.\end{Thm}
\begin{proof} We assume that $G$ doesn't have rank 1.
Let $I$ be a minimal closed invariant set for the action of $G$ on
$\partial X$. If
$radius _T(I)< \frac \pi 2$, then by \cite [II 2.7]{BRI-HAE}, $I$ has a unique centroid $b\in \bd X$.
Since $I$ is invariant,  $G$
fixes $b$,  and it follows by a result of Ruane \cite{RUA} (see
also Lemma \ref{L:point} of the present paper) that $G$ virtually
splits as $H\times \Bbb Z$. But then $diam_T(\bd X)= \pi $.  Thus we may assume that for
any $a \in \bd X$ and any $\delta >0$ there is a point of $y \in I$ with $d_T(a,y) > \frac \pi 2 -\delta$.

Let $x \in \bd X$.  Assume that $d_T(x,  I) =\frac \pi 2+\epsilon $ for
some $x\in \bd X$ and some $\epsilon >0$. Let $g_n\in G$ such that
$g_nt\to x$ for some (hence for all) $t\in X$. We may assume by
passing to a subsequence that $g_n^{-1}t\to n\in \bd X$. There is
some $y\in I$ such that $d_T(y,n)> \frac \pi 2 -\frac \epsilon 2$.

We apply now $\pi $-convergence to the sequence $g_n$ and the
point $y$. Clearly $y$ does not lie in the Tits ball of radius
$\frac \pi 2-\frac \epsilon 2$ and center $n$.
Passing to a subsequence,  we have that $g_ny\to z \in \bar B_T(x, \frac \pi 2+\frac \epsilon 2)$. Since $z \in I$ we have that $d_T(x,I)\leq
\frac \pi 2+\frac \epsilon 2$, a contradiction.

Thus we have shown that $d_T(x,I) \le \frac \pi 2$.  Choose $m\in I$ with $d_T(x,m) \le \frac \pi 2$.  Since $G$ doesn't have rank 1, by Theorem \ref{T:Bal}, $\bd X \subset \bar B_T(m, \pi)$.  By the triangle inequality, $\bd X \subset \bar B_T(x, \frac {3\pi} 2)$ as required.

\end{proof}

\begin{Cor} If the Tits diameter of $\partial X$ is bigger than
$\frac {3\pi}2 $, then $G$ contains a free subgroup of rank 2.
\end{Cor}

\section {The action on the cut-point tree is non trivial}
{\bf For the next two sections, we assume that $\bd X$ has a cut point.}
  We apply the construction of Section 1 to $Z= \bd X$ to obtain a non-trivial
$\mathbb{R}$-tree $T$ encoding the cut-points of $\bd X$. Since
the construction
  of $T$ is canonical, $G$ acts on $T$ by homeomorphisms.

   By the construction of $T$ if $G$ fixes a point of
 $T$ then $G$ fixes a point of $\cP$. Points of $\cP$ are of one
 of the following two types: i) cut points and ii)maximal inseparable sets
of $\partial
 X$.
  \begin{Lem}\label{L:nest} The action of $G$ on the $\mathbb{R}$-tree $T$
is non-nesting.  That is for any closed arc
   $I$ of $T$, there is no $g\in G$ with $g(I) \subsetneq I$.
  \end{Lem}
  \begin{proof}
   Suppose not, then by the Brouwer fixed point Theorem, there is an $A \in
T$ with $g(A)=A$.

  Replacing $g$ with $g^2$ if need be, we may assume that $I = [A,B]$ with
$g(B)  \in (A,B)$.
  We may also assume that  $A,B \in \cP$.  Similarly we may assume that
$|(g(B), B)|>2$, so $d_T(g(B),B) >0$.
   Since $G$ acts by isometries on the Tits boundary, it follows that
$d_T(B,A) = \iy$, moreover all limit points of the sequence of sets $\left (g^iB\right)$
are at infinite Tits distance from all limit points of the sequence of sets $\left (g^{-i}B\right)$.
So by $\pi$-convergence  $g$ is a rank 1 hyperbolic element.

   Now by \cite{BAL} every point of $\bd X$ is at infinite Tits distance
from both
   of the endpoints, $g^\pm \in \bd X$, of $g$.  By $\pi$-convergence
$g^+ \in A$. We consider $g^{-i}B$ and we note that as $i\to
\infty $ it converges either to a point $C$ of $T$ or to an end of
$T$ corresponding to an element $C$ of $\cP$. We remark that
$g^-$ lies in $C$ and that  $B \in (A, C)$.

Since $A$ is not terminal in $\cP$, there is a $D\in \cP$  with $A \in (D,C)$.  Thus  by
$\pi$-convergence  $g^{-i}D$ (for some $i>0$) lies in the
component of $T-A$ containing $C$. This is however impossible
since $A$ is fixed by $g$.

  \end{proof}

We deal now with the first type, i.e. we assume that $G$ fixes a
 cut point.
 We need the following unpublished result of Ruane.  We provide a
proof for completeness

\begin{Lem}\label{L:point} If $G$ virtually stabilizes a finite subset $A$ of $\bd X$, then
$G$ virtually has $\Z$ as a direct factor.
\end{Lem}
\begin{proof} Clearly there is a finite index subgroup $H<G$ which fixes $A$ pointwise.
Let $\{h_1, \dots h_n\}$ be a finite generating set of $H$. By
\cite{RUA}, the centralizer $Z_{h_i}$ acts geometrically on the
convex subset $Min(h_i)\subset Y$. By \cite{SWE} $A \subset
\Lambda Min(h_i) = \Lambda Z_{h_i}$ for all $i$.  Since $Z_{h_i}$
is convex, we have by \cite{SWE}  that the centralizer of $H$,
$Z_H= \cap Z_{h_i}$ is convex and
$$\bigcap \Lambda Z_{h_i} = \Lambda \left [\bigcap Z_{h_i} \right] \supset  A$$
  Thus $ A \subset  \Lambda Z_H$, so $Z_H$ is an infinite CAT(0) group, and  by
\cite{SWE} $Z_H$ has an element of infinite order.  Thus $H$
contains a central $\Z $ subgroup.  By \cite[II 6.12]{BRI-HAE},
$H$ virtually has $\Z$ as a direct factor.
\end{proof}

From this Lemma it follows that if $G$ fixes a point of $\bd X$
then $G$ is virtually $K\times \Z$. Thus $\bd X$ is a suspension
and $|\Lambda K|\ne 1$ (\cite{SWE}, Cor. p.345) so $\bd X$ has no
cut points, a contradiction.

We now prove the following:
\begin{Thm} \label{T:fmax}
$G$ fixes no maximal inseparable subset.
\end{Thm}
{\em Proof.} It follows from Theorem \ref{T:Ball}  \cite[Theorem A]{BAL-BRI} that if $G$ fixes a maximal
inseparable set then the Tits diameter of $\bd X$ is at most
$\frac 3 2 \pi$.

 Let's assume that $G$ fixes the maximal inseparable set $B$. We
 remark that if $x\in B$ then $Gx$ is contained in $B$ so the
 closure of $Gx$ is contained in $\bar B$. Thus
$\bar B$ contains a minimal invariant set that we denote by
 $I$.  By Theorem \ref{T:Ball}, $\bd X \subset \nb _T(\bar B, \frac \pi 2)$, the closed $\frac \pi 2$
Tits neighborhood of $\bar B$.

 If
$c$ is a cut point of $\bd X$ we have $\bd X-\{c\}=U\cup V$ with $U,V$
open and $B$ is contained, say, in $U$. With this notation we have
the following:
\begin{Lem} \label{L:circ2}There is a $2\pi$ geodesic circle
contained in $\bar V$.
\end{Lem}
\begin{proof}
Let $\alpha $ be a geodesic arc of length smaller than $\pi $
contained in $V$.
 Let $m$ be the midpoint of $\alpha $ and let $a,b$ be its endpoints.
Fix $t\in X$ and consider the geodesic ray $\gamma$ from $t$ to
$m$. Construct an increasing sequence $(n_i) \subset \N$ and
$(g_i) \in G$ such that
   $g_i(\gamma(n_i)) \to y \in X$, $g_i(m) \to m' \in \bd X$,
   $g_i(a) \to a' \in \bd X$ and $g_i(b) \to b'$.
    By \cite{BRI-HAE} $\angle_y(a',b') =
   d_T(a,b) \ge d_T(a',b') \ge \angle_y(a',b')$.  Thus
   $\angle_y(a',b')=d_T(a',b')$ which implies by the flat sector
   Theorem (\cite{BRI-HAE}, p.283, cor.9.9) that the sector bounded by the rays $[y,a')$ and $[y,b')$
   is flat. Let's denote the limit set of this sector by $\alpha
   '$.

    By passing to a subsequence if necessary we have that
   $g_i(y)$ converges to some $p\in \bd X$. Notice that by construction $g_i^{-1}(y) \to m$.  We remark that
   $d_T(p,m')\geq \pi$, so by $\pi $-convergence
   $g_i^{-1}(\alpha ')\cap V\ne \emptyset $ for $i>>0$. So there is a flat
   sector $Q$ (based at $v \in X$) whose limit set is contained in $V$.

   To simplify notation we denote the limit set of this sector
   again by $\alpha $ (and $a,b,m$ its endpoints and midpoint respectively).

Repeating this process (by pulling back along the ray
   $[v,m)$) we obtain  a sequence $(h_i) \subset G$ and a  flat $F$ with $h_i(Q) \to F$ (uniformly on compact subsets of $F$), and so a geodesic $2 \pi$-circle $\gamma= \Lambda F$ with
   $h_i(\alpha)$ converging to an arc of $\gamma$.
   As before, using $\pi $-convergence, we show that for $i>>0$, $h_i^{-1}(\gamma) \cap V\neq \emptyset$ and since $\gamma$ is a {\bf simple} closed curve, $h_i^{-1}(\gamma) \subset \bar V$.
\end{proof}

Since no geodesic circle can lie in a Tits ball of radius less than $\frac \pi 2$, and $X \subset \nb _T(B, \frac \pi 2)$, it follows from Lemma \ref{L:circ2} that:
\begin{itemize}
\item $c \in \bar B$
\item $\exists p\in V$ such that $d_T(p,\bar B)=\pi /2$.
\end{itemize}
Notice that we have shown that every cut point of $X$ is adjacent to $B$ (in the pre-tree structure) and so every  cut point is contained in $\bar B$.

 Fix $t\in
X$ and consider a sequence $(g_i) \in G$ such that
  $g_i(t) \to p \in \bd X$. By passing to a subsequence we may
  assume that $g_i^{-1}(t) \to n \in \bd X$. Since $G$ does not
  virtually fix $c$, there are distinct translates $h_1(c),h_2(c) \neq n$
which do not separate $B$ from $n$.  Thus
  $d_T(h_1(p), n) > \frac \pi 2$ and $d_T(h_2(p), n) > \frac \pi 2$. By $\pi $-convergence $g_i(h_1(p)),
  g_i(h_2(p))$ lie in $V$ for sufficiently large $i$. Since all cut points are adjacent to $E$, we have $g_i(h_1(c))=g_i(h_2(c))=c$ for sufficiently
  big $i$. This is impossible as $h_1c,h_2c$ are distinct, and the proof of Theorem \ref{T:fmax} is complete.\qed

 \begin{Cor}\label{C:rank1} In our setting, the action of $G$ on $X$ is rank 1, so $\bd X$ has infinite Tits diameter.
 \end{Cor}
 \begin{proof} The action of $G$ on the $\R$-tree $T$ is non-nesting and  without global fixed points.
 It follows from \cite[Proposition 35]{P-S} that $G$ has an element $h$ which acts by translation on a line of $T$.
 Thus there is a cut point $c$ such that $c$ separates $h^{-1}(c)$ from $h(c)$.  It follows by $\pi$-convergence that $d_T(h^+,h^-) = \iy$.
 \end{proof}

\section {Cat(0) groups have no cut points}

We study now further the action of $G$ on $T$.
\begin{Lem} \label{L:finitestab}
If an interval $I \subset \cP$ is infinite then the stablizer of $I$ in $G$ is finite.
\end{Lem}
 \begin{proof}
Assume that $I$ is fixed (pointwise) by an infinite subgroup of
$G$, say $H$. Let $(h_i)\subset H$ be an infinite sequence of
distinct elements of $H$. By passing to a subsequence we may
assume that $h_i(t)\to p\in \bd X$, $h_i^{-1}(t)\to n\in \bd X$
for any $t \in X$.

The interval $I$ contains infinitely many cut points, so by Corollary \ref{C:rank1} we can find $x \in \bd X$,  $d_T(n,x) = \iy$,
  with the property that $x$ is separated from $p$ by two cut points $c,d \in I$.  We may assume that
  $d$ separates $c$ from $p$.
 By $\pi $-convergence  $h_i(x)\to p$.  There exist subcontinua $B \ni p$ and $A \ni c,x$ with $X = A \cup B$ and
 $A \cap B =d$.  Furthermore, since $(h_i)$ fixes $c$ which separates  $x$ from $d$, it follows
 that $c$ separates $h_i(x)$ from $d$.   Thus $\{h_i(x) \} \subset A$,  contradicting $h_i(x) \to p \in B - \{d\}$.

\end{proof}

\begin{Lem}\label{L:stab} The action of $G$ on the $\mathbb{R}$-tree $T$
is stable.
  \end{Lem}

  \begin{proof}
We recall that a non-degenerate arc $I$ is called stable if there
is a (non degenerate) arc $J\subset I$ such that for any non
degenerate arc $K\subset J$, $stab(K)=stab(J)$.

An action is called stable if any closed arc $I$ of $T$ is stable.

 We remark that by construction, arcs of $T$ that correspond
to adjacent elements of $\cP $ are stable. Further if an arc
contains a stable arc it is itself stable, so every unstable arc
contains infinitely many elements of $\cP$.

Suppose now that an arc $I$ is not stable.   Then there is a
properly decreasing sequence $I=I_1\supset I_2\supset ... $ so
that $stab(I_1)\subset stab(I_2)\subset ...$ where all inclusions
are proper.   Notice that $I_i \cap \cP$ is infinite for each $i$,
and so $stab(I_i)$ is finite by Lemma \ref {L:finitestab}. On the
other hand, there is a uniform bound on the size of a finite
subgroup of $G$, which is a contradiction.

\end{proof}

\begin{Lem}\label{L:discr}  $T$
is discrete.
  \end{Lem}

  \begin{proof}
  Suppose that $T$ is not discrete.  Then by \cite{Le} there is an
$\R$-tree $S$
 and a $G$-invariant  quotient map $f:T \to S$  such that $G$ acts
non-trivially  on $S$ by isometries.
 Furthermore for each non-singleton arc $\alpha \subset S$,
$f^{-1}(\alpha)\cap \cP$ is an infinite interval
 of $\cP$, and the stabilizer of an arc of $S$ is the stabilizer of an
arc of $T$.
 It follows that the arc stabilizers of $S$ are finite (Lemma \ref{L:finitestab}) and that the action
of $G$ on $S$ is stable (Lemma \ref{L:stab}).

 Since $G$ is one-ended, by the Rips machine applied to $S$,
 there is an element $h\in G$ which acts by translation on a line in $S$
such that $G$ virtually splits over $h$.
 This implies that $\{h^\pm\}$ is a cut pair of $\bd X$, however since $h$ acts by
translation on a line in $S$,
 it acts by translation on a line in $T$, and so $h^\pm$ will correspond
to the ends of this line (by $\pi$-convergence).
 Thus $\{ h^\pm\}$ lie in terminal elements
 of $\cP$.  This contradicts the fact that $\{h^\pm\}$ separates.
 Thus $T$ is discrete.

  \end{proof}

 The following is easy and well known, but we provide a
proof for completeness
\begin{Lem}\label{L:centlim}  If $H<G$ then the centralizer $Z_H$ fixes
$\Lambda H $ in $\bd X$.
\end{Lem}
\begin{proof}
Let $\alpha \in \Lambda H \subset \bd X$. By definition there
exists a sequence of elements $(h_n)\subset H$ with $h_n(x) \to
\alpha$ for any $x \in X$.  Since $G$ acts by homeomorphisms on
$\bar X = X \cup \bd X$, $g(h_n(x)) \to g(\alpha)$ for any $g\in
Z_H$. Notice however that $g(h_n(x)) =h_n(g(x))$, and $h_n(g(x))
\to \alpha$.  It follows that $\alpha = g(\alpha)$.
\end{proof}

For a finitely generated $H <G$, the centralizer $Z_H$ will be a {\em convex} subgroup of $G$, that is $Z_H$ will act geometrically on a convex subset of $X$.
\begin{Lem} \label{L:fixadj} Let $H<G$ with $H$ fixing the adjacent pair
$\{c,B\}$ of  $\cP$ (so $c$ is a cut point,
 $c \in \bar B$, $h(c)=c$ and $h(B) =B$ for each $h \in H$).  Then
$\Lambda H \subset \bar B$.
\end{Lem}
\begin{proof}   Let $(h_i) \subset H$ be a sequence with $h_i(x) \to p \in
\Lambda H - \bar B$ for any $x \in X$.
   Passing to a subsequence we have $h^{-1}_i(x) \to n \in \Lambda H$.

  The adjacent pair $\{c,B\}$ gives a separation of $\bd X- \bar B$ in the obvious
way namely:
   Let $U = \{\beta \in \bd X - \bar B: c \in ([\beta], B)\} $ and
$V=\{\beta \in \bd X - \bar B: B \in ([\beta], c)\} $.
     Using the definition of $\cP$  one can show  that both $U$ and $V$
are open in $\bd X$ (being the union of open sets).
  Using the pretree axioms we see that $U \cap V = \emptyset$ and $U \cup
V = \bd X-\bar B$.  Since $H$ leaves  $B$ and $c$
  invariant, $h_i(U) =U$ and $h_i(V) = V$ for all $i$.

 Without loss of generality let  $p \in U$.
Since  $G$ is rank 1 in its action on $X$, there is a point
$\alpha \in V$ with   $d_T(\alpha, n) >\pi$.  It follows by
$\pi$-convergence that $h_i(\alpha) \to p$. This contradicts the
fact that $p \in U$ and $h_i(\alpha) \in V$ (in particular
$h_i(\alpha) \not \in U$) for all $i$.  Thus $\Lambda H \subset
\bar B$.
  \end{proof}
\begin{Lem}\label{L:adinf}
Suppose that we have  adjacent pairs $\{B,c\}$ and $\{\hat B,c\}$
in $\cP$ with $B \neq \hat B$ with the stabilizer  of $\{B,c\}$
infinite and the stabilizer of $\{\hat B,c\}$  infinite.  Then
there is a hyperbolic $g \in G$ with $g(B) =B$, $g(\hat B) = \hat
B$, and $g(c) =c$.
\end{Lem}
\begin{proof}
Since there are only finitely many conjugacy classes of finite
subgroups in $G$, there are infinite finitely generated subgroups
$H= \langle h_1\dots h_n\rangle$ and   $\hat H = \langle \hat
h_1\dots \hat h_m\rangle$ stabilizing $\{B,c\}$ and  $\{\hat B,
c\}$ respectively. Since $h_i(c) =c =\hat h_j(c)$ for all $i,j$ it
follows from the \cite{SWE} that $c \in \Lambda Z_{h_i} $ and $c
\in \Lambda Z_{\hat h_j}$ for all $i,j$. Since $Z_{h_i}$ and
$Z_{\hat h_j}$ are convex, it follows from \cite[Theorem 16]{SWE}
that $Z_H = \bigcap\limits_i Z_{h_i}$ and $Z_{\hat H}
=\bigcap\limits_i Z_{h_i}$ are convex and that $c \in \Lambda Z_H$
and $c \in \Lambda Z_{\hat H}$.  Thus by \cite[Theorem 16]{SWE} $Z
= Z_H \cap Z_{\hat H}$ is convex and $c \in \Lambda Z = \Lambda
Z_H \cap \Lambda Z_{\hat H}$.

Since $Z$ is convex, it is a CAT(0) group, and by \cite[Theorem
11]{SWE} there is an element $g \in Z$ of infinite order.
 For any $K <G$, $\Lambda K$ is
not a single point \cite[Corollary to Theorem 17]{SWE}. Since $H$ is infinite, $\Lambda H$ is
non-empty with more than one point.  By Lemma \ref{L:fixadj}
$\Lambda H \subset \bar B$.  By Lemma \ref{L:centlim},  $g$ fixes
$\Lambda H$.  It follows that $g(B) =B$.  Similarly $g(\hat B)
=\hat B$.  It follows that $g(c) = c$ as required.
   \end{proof}

\begin{Thm} \label{T:nocut} Let $G$ be a one-ended group acting
geometrically on the CAT(0) space $X$. Then $\bd X$ has no cut
points.
\end{Thm}
\begin{proof} Assume that $\bd X$ has a cut point. Let $T$ be the
cut point tree of $\bd X$. By Lemma \ref{L:discr} $T$ is a simplicial
tree. Since $G$ is one-ended, all edge stabilizers are infinite. By
Lemma \ref{L:adinf} there is a hyperbolic element $g\in G$
  stabilizing two adjacent edges of $\{B,c\}$ and $\{c,D\}$ of $T$. Thus by  Lemma \ref{L:fixadj}, $\{g^\pm \} \subset \bar B$ and $\{g^\pm\} \subset \bar D$.  This is absurd since the intersection of $\bar B$ and $ \bar D $ is a single point (namely $c$).
\end{proof}

\section{The action on the JSJ-tree is non-trivial}
\begin{Se} {\bf For the remainder of the paper $G$ will be a one-ended group
  acting geometrically on the CAT(0) space $X$.  As we showed in the
previous section, $\bd X$
has no cut points.  We will further assume that $\bd X$ contains a cut pair.}
\end{Se}

Let $T$ be the JSJ-tree of $Z= \bd X$ constructed in section 1.
   Since the construction
  of $T$ is canonical, $G$ acts on $T$.
  \begin{Lem}\label{L:nesti} The action of $G$ on the $\mathbb{R}$-tree $T$
is non-nesting.  That is for any closed arc
   $I$ of $T$, there is no $g\in G$ with $g(I) \subsetneq I$.
  \end{Lem}
  \begin{proof} The proof is basically the same as that of Lemma \ref {L:nest}.  We sketch it here.
   Suppose not, then we may assume that:
   \begin{itemize}
   \item $I = [A,B]$
   \item $A,B \in \cP$
   \item  $g(A)=A$
   \item $g(B)  \in (A,B)$
\item $|(g(B), B)|>2$, so $d_T(g(B),B) >0$
\end{itemize}
Using $\pi$-convergence we see that :
\begin{itemize}
\item $d_T(B,A) = \iy$.
\item $g$ is a rank 1 hyperbolic element.
 \item
$g^+ \in A$.
\item $B$  separates $A$ from
$g^-$
\end{itemize}
No element of $\cP$ is a singleton, so let $a \in A -\{g^+\}$.
Notice that $g$ is rank 1 so $d_T(a, g^+) = \iy$.
   Thus by $\pi$-convergence $g^{-i}(a) \to g^-$.
    Since $A$ is fixed by $g$,  then $g^- \in \bar A$.  This contradicts  the fact that $B$  separates $A$ from
$g^-$ \end{proof}

 By the construction of $T$ if $G$ fixes a point of
 $T$ then $G$ fixes a point of $\cP$. Points of $\cP$ are of one
 of the following types:
 \begin{itemize}
 \item inseparable cut pairs of $\partial X$
 \item maximal
 cyclic subsets of $\partial X$ which we call {\em necklaces}
\item  maximal inseparable subsets
 of $\partial X$.
 \end{itemize}

 We deal now with the first type, i.e. we assume that $G$ leaves an
 inseparable cut pair invariant, so in fact $G$ virtually fixes each element of the cut pair.
 By Lemma \ref{L:point}
 if $G$ virtually fixes a point of  $\bd X $  then
$G$ is virtually $ H \times \Z$ for some finitely presented group.  Since $\bd X$ has a cut pair, $H$ will have more than one end.  In this case either $G$ is virtually $\Z^2$ or $\bd X$ has only one cut pair and all the other elements of $\cP$ are maximal inseparable sets.  Having classified these possibilities, we now ignore them.

\begin{Se} {\bf For the remainder of the paper, we assume that $G$ doesn't
virtually have $\Z$ as a direct factor,  in particular $G$ doesn't leave  a cut pair (or any other finite subset) invariant}.
\end{Se}

\subsection{Maximal inseparable subsets of $\partial X$
}

\begin{Lem} If the Tits diameter of $\bd X$ is more than $\frac {3\pi} 2$,
then $G$ leaves no  maximal inseparable set invariant.
\end{Lem}
\begin{proof}  If $A $ is a cut pair, then there are nonempty disjoint
open subsets $U$ and $V$ such that $U\cup V = \bd X-A$.  Choose a
rank 1 hyperbolic element $g \in G$ with $g^- \in U$ and $g^+ \in
V$.   By $\pi$-convergence, it follows that $g(A) \subset V$ and
$g(A)$ separates $A$ from $g^+$.   By $\pi$-convergence every
maximal inseparable set will lie "between" $g^i(A)$ and
$g^{i+1}(A) $ for some $i \in Z$, and so will  be moved off itself
by $g$.
\end{proof}

Let's assume now that $G$ leaves invariant a maximal inseparable subset of
$\partial X$, say $D$, so the Tits diameter is at most $\frac
{3\pi}2$. Let $a,b$ be a pair of points separating $\partial X$.

Let $B$ be the closure of the component of $\partial X-\{a,b\}$
which contains $D$, and $I$ be a minimal non-empty closed
invariant set for the action of $G$ of necessity contained in $D$.
By \cite[Proposition II 2.7]{BRI-HAE} if the  Tits radius of $I$
were less than $\frac \pi 2$, then $I$ would have a unique centroid (in
$\bd X$) which would be fixed by $G$.  Since $G$ fixes no point of
$\bd X$, it follows that the Tits radius of $I$ is at least $\frac
\pi 2$. Similarly $I$ must be infinite.

\subsubsection{$C$-geodesics} In order to show that $G$ leaves no maximal inseparable set invariant, we must understand a special type of local Tits geodesic in $\bd X$.
Let $\{c,d\}\in \partial X$ be a cut pair and    $C \subset \bd X$
the closure of a component of  $\bd X -\{c,d\}$. Consider the path
metric of $C$ using the Tits metric $d_T^C(x,y) = \inf \{
\ell_T(\alpha) : \alpha$ is a path in $C$ from $x$ to $y \}$.  For
$e,f \in C$, let $\epsilon >0$, consider the sets $F_\epsilon =
\cup \{ $ paths $\alpha \subset C$ with $f \in \alpha $ and with
Tits length $\ell_T(\alpha ) \le \epsilon \}$ and $E_\epsilon =
\cup \{ $ paths $\alpha \subset C$ with $e \in \alpha $ and with
Tits length $\ell_T(\alpha ) \le \epsilon\}$. Using lower
semi-continuity of the (identity) function from the Tits boundary
$TX$ to the boundary $\bd X$, we see that $F_\epsilon$ and
$E_\epsilon$ are closed connected subsets of the continua $C$.  The Tits
diameter of $\bd X $ is at most $\frac {3 \pi} 2$.  It follows that for some
$\epsilon $, $ C= E_\epsilon \cup F_\epsilon$, so
 there exists a path $\alpha \subset C$
of finite Tits length from $f$ to $e$. Using lower semi-continuity
of the identity map (from the Tits boundary to the regular
boundary) and a limiting argument, we can show that there is a
shortest path in $C$ from $f$ to $e$ which we call a $C$-geodesic
from $f$ to $e$.

Local geodesics of length $\le \pi$ are geodesics, \cite{BRI-HAE},
so if $d_T^C(d,e) <\pi$, then the $C$-geodesic from $d$ to $e$
is the Tits geodesic from $d$ to $e$ and therefore unique.

\begin{Thm}
$G$ doesn't leave a maximal inseparable set invariant.
\end{Thm}
\begin{proof}

Let's assume that $G$ leaves invariant a maximal inseparable subset of
$\partial X$, say $D$. Let $I$ be a minimal invariant set for the
action of $G$ on $\bd X$.  Observe that $I$ is contained in the closure of $D$,
 $\bar D$.
It follows that $\bd X$ has finite Tits diameter.

$D$ is a point of $T$. Let $R$ be a component of
$T-D$.  We can view $R$ as a subset of $\bd X$ in the appropriate way.  With this viewpoint,  the closure of $R$ in the Tits topology (i.e. the topology defined by $d_T$)
intersects $D$ at either one or two points.

We distinguish two cases:
\begin{enumerate}
\item There is no element of $\cP$ contained in $R$ adjacent to $D$.
\item There is an element of $\cP$ contained in $R$ adjacent to $D$.
\end{enumerate}
In the first case
there is a non trivial loop  contained in $R$  separated from $D$ by infinitely many elements of $\cP$.
 We homotope this loop  to
a geodesic circle $w$  separated from $D$ by infinitely many elements of $\cP$.  Thus we have  $\{\hat a,\hat b\}$, a cut pair disjoint from $w$ and $D$,
 separating $w$ from $D$.

 Then we claim that for
some $e\in w$, $d_T(e, \{\hat a, \hat b\}) \ge \frac \pi 2$.
Indeed if every point on $w$ is at Tits distance less than $\frac
\pi 2$ from $\hat a, \hat b $ then there are two antipodal points
on $w$, both at Tits distance less than $\frac \pi 2$ from one of $\hat
a, \hat b $. This contradicts the fact that $w$ is geodesic.
It follows now that $d_T(e,I)> \frac \pi 2$ which is a
contradiction.

We deal now with the second case. Since the Tits diameter of $\bd
X$ is at most $\frac {3\pi} 2$, it follows that the element of
$\cP$ adjacent to $D$ in $R$ is a cut pair, $\{a,b\}$.  All
translates (there are infinitely many of them) of $\{a,b\}$ are
adjacent to $D$. Let $p$ be a point of $\bd X$ separated from $D$
by $\{a,b\}$,  and let $g_i\in G$ such that for some (hence for
all) $x\in X$, $g_i(x)$ converges to $p$. By passing to a
subsequence we may assume that $g_i^{-1}x$ converges to some $n\in
\partial X$.

By \cite{ONT}, $X$ has almost extendable geodesics, and this
implies that there is $q \in \bd X$ with
 $ d_T(n,q)=\pi$,
 and by $\pi $-convergence $g_i(q) \to p$.
It follows that there are translates $a_1,b_1$ of
 $a,b$ respectively  so that  $q $ is separated from $D$ by $\{a_1,b_1\}$.
Further we claim that $g_i\{a_1,b_1 \}=\{a,b\}$ for all $i$ big
enough. Indeed since $\{a,b \}$ is a cut pair there is a
neighborhood $U$ of $p$ with $\bd U = \{a,b\}$. On the other hand
if $g_i\{a_1,b_1 \}\ne \{a,b\}$, $g_i(q) \not \in U$.
 We have in fact shown that  every point of $\bd X$ at Tits distance $\ge
\pi$ from $n$ is separated from $D$ by $\{a_1,b_1\}$.

 Let $B_1$ be the component of $\bd X-\{a_1,b_1\}$ containing $D$ and
$C_1$ the component of
 $\bd X-\{a_1,b_1\}$ containing $q$. We claim that $n\in B_1$.
 Indeed suppose that $n$ lies in another component, say $B_2$ of $\bd
 X-\{a_1,b_1\}$. By Theorem \ref{T:Ball}, every point is at Tits distance less than
 $\frac \pi 2$ from $I$, so there is a $B_2$-geodesic $\gamma _1$ of length less than
 $\pi $ from $a_1$ to $b_1$.  It follows that if $\gamma _2$ is a
 $B_1$ geodesic from $a_1$ to $b_1$ the length of $\gamma _2$ is
 at least $\pi $. This is because $\gamma _1 \cup \gamma _2$ is a
 non-contractible loop. It follows that there is some point on $\gamma _2$
 at distance more than $\pi $ from $n$. This contradicts our
 earlier observation that all such points are separated from $D$
 by $\{a_1,b_1\}$.

 Now consider the loop $S$ consisting of the $B_1$-segments $[n,a_1]$,
$[n,b_1]$ and a
 $C_1$-segment $[a_1,b_1]$. We may assume that $q\in [a_1,b_1]$ because $S$ is
non-contractible, so some point of $S$ is at Tits distance $\ge \pi$ from $n$.
 The $B_1$-segments $[n,a_1]$ and $[n,b_1]$, having length less than $\pi$,
will be Tits geodesics.   Furthermore, since every point of $\bd X$ is within $\frac \pi 2$ of $I \subset D$,
 the $C_1$-segment $[a_1,b_1]$ must also have length less than $\pi$ and is therefore a Tits geodesic.

 Since every essential loops of $\bd X$ has Tits length at least $2\pi$, one can show that only finitely many translates of $\{a,b\}$ separate a subarc of $S$ from $D$.
Let $\{a_2,b_2 \}$, $C_2,B_2$ be translates of $\{a_1,b_1 \}$,
$C_1,B_1$ respectively so that $\{a_2,b_2\}$ doesn't separate a
subarc of $S$ from $D$.  As before let $S'$ be the union of the
$B_2$-segments $[n,a_2]$,  $[n,b_2]$ and a
 $C_2$-segment $[a_2,b_2]$, and as before these three segments will be
Tits geodesics.  Notice that since $S'$ is
 not contractible, there is a point of $S'$ at distance $\ge \pi$ from
$n$,  and so some point of $S'$ is separated from $D$ by $\{a_1,b_1\}$

Notice by construction that for every $s \in S$, $S$ contains a
Tits geodesic from $s$ to $n$ and similarly for $S'$. Then $S,S'$
both contain $a_1,b_1$, have length at least $2\pi $. Let
$S_1,S_1'$ be respectively the subarcs of
 $S,S'$ from $n$ to $a_1$ contained in $B_1$ and similarly let $S_2,S_2'$
be respectively the subarcs of
 $S,S'$ from $n$ to $b_1$ contained in $B_1$.
Then as we note before all have length smaller than $\pi $. On the
other hand at least one of $S_1\cup S_1'$, $S_2\cup S_2'$ runs
through $C_2$ once, and so is a non contractible loop. This is a
contradiction.

\end{proof}

\subsection{Stabilizers of  necklaces of $\partial X$}

\begin{Thm} \label{T:onto} Let $H$ be a group of homeomorphisms of the
circle $S^1$, and let $A$ be a nonempty  closed $H$-invariant
subset of $S^1$.
 Either $A$ is countable and $H$ virtually fixes a point of $A$,
 or there is an $H$-equivariant cellular quotient map $q: S^1 \to S^1$ such
that $q(A) = S^1$.
\end{Thm}
\begin{proof}
Recall that  a compact Hausdorff space $E$ is called perfect if
each point of $E$ is a limit point of $E$. We note that every
perfect space is uncountable.

 Let  $$B = \{ a \in A : \exists \text{ an open  }U\ni a \text{ with $U \cap
A$ countable }\}$$
 and let $C = A-B$.  Since $A$ is Lindeloef, and $B$ is a subset of $A$, $B$
is countable. Since $B$ is an open subset of $A$, $C$ is compact,
and it follows that $C$ is perfect since every neighborhood of a
point of $C$ contains uncountably many points of $A$, and
therefore  a point of $C$ (since $B$ is countable).  Thus $A$ is
the union of a countable set $B$, and a perfect set $C$.

First consider the case where $A$ is countable, ($C= \emptyset$).
We remark that if $A\supset A_1\supset A_2\supset ...$ and $A_i$
are closed $H$-invariant sets, then $\cap A_i$ is non-empty and
$H$-invariant. By Zorn's Lemma there is a minimal, closed,
non-empty, $H$-invariant subset of $A$. Let's call this set $A'$.
If $a\in A'$ then $A'$ is the closure of the orbit of $a$ under
$H$, $\overline {Ha}$. If $A'$ is infinite then there is some
limit point $b$ in $A'$. Since $A'$ is minimal closed
$H$-invariant set, and $\overline {Hb}\subset A'$ we have that
$\overline {Hb}= A'$. It follows that $a$ is a limit point of
$A'$. Since $A'=\overline {Ha}$, every point of $A'$ is a limit
point of $A'$, so $A'$ is perfect, a contradiction. Therefore $A'$
is finite and $H$ virtually fixes a point of $A$.

Now consider the case where $A$ is uncountable, so $C \neq
\emptyset$.

 Let $I$ be a closed interval of $\R$ and let $D$ be a perfect subset of
$I$.  For $x,y \in I$ define $x\sim y$ if there is no point of $d$
strictly between them.  This is an equivalence relation
(transitivity follows from the perfection of $D$).
 The linear order of $I$ extends to a linear order on the quotient $I/\sim$.
It follows that $I/\sim$ is a compact  connected separable
linearly ordered topological space, also known as a closed interval.
Since the  image of $D$ in $I/\sim$ is compact and dense, it
follows that the image is the entire quotient $I/\sim$.

 Of course $S^1$ is just the quotient space of a closed interval which
identifies exactly the endpoints.
 We now say that $x,y \in S^1$ are equivalent if there is an open arc in
$S^1-C$ from $x$ to $y$.
 It follows from the previous paragraph that this is an equivalence relation
on $S^1$ and that $S^1/\sim = S^1$.  Let $q:S^1 \to S^1$ be the
quotient map.  Since $C$ is $H$-invariant, it follows that  $q$ is
$H$-equivariant, and by the previous paragraph $q(C) = S^1$.  By
constuction, the inverse image of a point is a closed  arc, so $q$
is cellular.
\end{proof}

\begin{Lem}\label{L:rankone} If $G$ leaves a necklace $N$ of $\bd X$ invariant, then the action of $G$ on $X$ is rank 1.
\end{Lem}
\begin{proof}  Suppose not. Since the Tits diameter of $\bd X$ is at most $\frac 3 2 \pi$, $N$ is contained in a geodesic circle $\alpha$ of length at most $ 3  \pi$.  If any arc of $\alpha-N$ has length at least $\pi$, then the endpoints of that arc are virtually stabilized by $G$ which contradicts our hypothesis.  Thus  $\alpha$ is the unique geodesic circle containing $N$, and so  $G$ stabilizes $\alpha$.

Since $G$ acts on the Tits
boundary by isometries, there is a homomorphism $\rho:G \to
Isom(\alpha)$ (which is virtually the Lie group $S^1$ and so virtually abelian).  Since $G$ doesn't virtually fix any point of $\bd X$,  $\rho(G)$ is infinite.  Since
$G$ is finitely generated, the virtually abelian group $\rho(G)$ contains an element of
infinite order, $\rho(g)$.  That is $g$(or $g^2$) rotates $\alpha$ by
an irrational multiple of $\ell_T(\alpha)$.  It follows by $\pi$-convergence that $d_T(s, g^+) =d_T(s,g^-) =\frac \pi 2$ for all $s \in \alpha$.  Thus every point of $N$ is joined to $g^+$ by a Tits geodesic of length exactly $\frac \pi 2$. However since $N$ is a necklace and $g^+ \not \in N$,  there must be $a,b ,c \in N$ such that $\{a,b\}$ separates $c$ from $g^+$, and so $d_T(c,g^+) > \min\{d_T(a, g^+), d_T(b,g^+)\}$, a contradiction.
\end{proof}

\begin{Cor}\label{C:rank11} In our setting ($\bd X$ has a cut pair and $G$ doesn't virtually have $\Z$ as a direct factor) the action of $G$ on $X$ is rank 1.
\end{Cor}
\begin{proof}
By Lemma \ref{L:rankone} we may assume that $G$ doesn't leave a necklace of $\bd X$ invariant.
Thus $G$ acts on $\cP$ without fixing a point, so $T\neq \emptyset$ and $G$ acts on $T$ without fixed points.  By  \cite{P-S}
that there is $h \in G$ which acts as a hyperbolic element on $T$, that is
$h$ acts by translation on a line $L \subset T$. Let $A \in L \cap
\cP$.  For some $n>0$, $d_T(h^n(A) , A) >0$.  It follows by
$\pi$-convergence that $d_T(h^+, h^-) = \iy$, and by definition
$h$ is rank 1 in its action on $X$.
\end{proof}
\begin{Def} Let $S$ be a necklace of $X$.  We say
$y,z \in X-S$ are $S$ equivalent, denoted $y\sim_S z$ if for every cyclic
decomposition $M_1, \dots M_n$ of $X$ by $\{x_1, \dots x_n \} \subset S$,
both $y,z \in M_i$ for some $1\le i \le n$.  The relation $\sim_S$ is
 an equivalence relation on $X-S$.
The closure (in $X$) of a $\sim_S$-equivalence class of $X-S$ is called a
{\em gap} of $S$.  Notice that every gap is a nested intersection of
continua, and so is a continuum.
Every inseparable cut pair in $S$ defines a
unique gap. The converse is true if $X$ is locally connected, but false in
the non-locally connected case.
\end{Def}
\begin{Thm} \label{T:necklace} Let $B$ be a necklace of $\bd X$, and
 let $H$ be the stabilizer of $B$, so $H = \{g \in G: g(B) = B\}$.
 Then one of the following is true:
\begin{enumerate}
\item $H$ is virtually cyclic and there is at most one gap $J$ of $B$ with
the stabilizer of $J$ in $H$ infinite .
\item $H$ acts
properly by isometries on $\HH^2$ with limit set the circle at
infinity. Furthermore for any gap $J$ of $B$, the stabilizer of
$J$ in $H$ is finite or it is a peripheral subgroup of $H$ (as a
Fuchsian group). Distinct gaps correspond to distinct peripheral subgroups.
\end{enumerate}
\end{Thm}
\begin{proof}

  Let $f_B:\bd X \to S^1$
be the circle function for $B$ (see Theorem 22 of \cite{P-S}).
Thus the action of $H$ on $B$ extends to an action of $H$ on
$S^1$.  By Theorem \ref{T:onto} either $H$ virtually fixes a point
in $\overline {f_B(B)}$, or by composing $f_B$ with an
$H$-invariant cellular quotient map $q$, we have an $H$ invariant
$\pi: \bd X \to S^1$. We may assume that $H$ is infinite.

\subsubsection*{Case I}
We first deal with the case where $H$ virtually fixes a point $r
\in  \overline {f_B(B)}$. Passing to a finite index subgroup, we
assume $H$ fixes $r$.

First consider the case where an infinite subgroup $K < H$ fixes three points
$a,b,c \in B$.  We have the  unique cyclic decomposition
$U,V,W$, continua with $U\cup V \cup W = \bd X$ and  $U\cap V= a$, $V \cap W =b$, and $W\cap U = c$.   With no loss of generality, there is $p \in \Lambda K -U$.  Choose $(k_i) \subset K$ with $k_i(x) \to p$ and $k_i^{-1}(x) \to n$ for some $n \in \bd X$ and for all $x \in X$.   Since $G$ is rank one, there exists
$q \in U$ with $d_T(q,n) =\iy$.  Thus $k_i(q) \to p \not \in U$, but since $K$ fixes $a,b,c$, by uniqueness of the cyclic decomposition $K$ leaves  $U$ invariant, a contradiction.  {\em Thus no infinite subgroup of $H$ fixes three points of $B$; in particular no hyperbolic element of $H$ fixes more than two points of $B$}.

We remark now that $ \overline {f_B(B)}-p$ is linearly ordered so
a finite order element of $H$ either fixes all points of $
\overline {f_B(B)}$ or its square fixes all points of $ \overline
{f_B(B)}$.  If  $H$ is torsion, then passing to an index two subgroup, we may assume that
$H$ fixes $ \overline{f_B(B)}$, so $H$ fixes $B$.  As noted above, this is not possible, so $H$ contains a hyperbolic element.

 Since $H$ fixes a point in $S^1$, the action of $H$ on $S^1$ comes from an
action on a closed interval.
 Thus there exists $b \in B$  and hyperbolic $h \in H$ such that $\{b,
h^2(b)\}$ separates $h^{-1}(b) $ from $h(b)$.

 Since $d_T(b,g(b))>0$,
for each $i$, $h^i(b)$ is at infinite Tits distance from any fixed
point of $h$ in $\bd X$, in particular from $h^{\pm}$. It follows
by $\pi$-convergence that $h^i(b) \to h^+$ and $h^{-i}(b) \to
h^-$. For $m\in \Z$ let $M_m$ be equal to the closure of the
component of $\bd X - \{h^{m-1}(b), h^{m}(b)\}$ which doesn't
contain $h^{\pm }$.  Notice that $h(M_m) =M_{m+1}$, and that
$d_T(M_m, h^{\pm })=\iy$. From $\pi$-convergence, the closure of
$M=\cup M_m$ is $M \cup \{h^\pm\}$.  It follows by \cite[Lemma
15]{P-S} that $B \cup \{h^\pm\}$ is cyclic subset of $\bd X$, and
by maximality, $\{h^\pm\} \subset B$. Since $h$ fixes at most 2
points of $B$, $h^\pm$ are  the only points of  $B$ fixed by $h$,
and  $r \in  f_B(h^\pm)$.  Since $f^{-1}_B( f_B(h^\pm))=\{h^\pm\}$, it
follows that $H$ fixes one of $h^\pm$. From \cite{BRI-HAE}, $H$
fixes both $h^\pm$.

We now show that $\Lambda H = \{h^\pm\}$.  Suppose not, then there
is $(g_i) \subset H$ with $g_i(x) \to p \not \in \{h^\pm\}$ and
$g_i^{-1}(x) \to n $ for all $x \in X$.  Notice that $(g_i)$ leaves $M$ invariant.  Since
$(g_i)$ fixes $h^\pm$, by $\pi$-convergence,
 $n,p\in B_T(h^\pm, \pi)$.   It follows that $d_T(M_m, n) =\iy$ for all
$m \in \Z$.  Thus $g_i(M_0) \to p$.  However $g_i(M_0) \subset M$
for all $i$ and the closure of $M$ is $M \cup \{h^\pm\}$ a
contradiction.  We have shown that  $\Lambda H = \{h^\pm\}$ and by  Lemma
\ref{L:limsets}, $H$ is virtually $\langle h \rangle$.  If there
is a gap of $B$ with sides $h^\pm$ then that gap will be
stabilized by $\langle h \rangle$
 (there can be at most one such gap), but no other gaps of $B$
 will be stabilized by a positive power of $ h$, and so will have finite
stabilizer in $H$.

\subsubsection*{Case II} By Theorem \ref{T:onto} we are left with the case where there is an $H$
invariant function $\rho: \bd X \to S^1$ with $\overline{\rho( B)}
= S^1$, where $\rho = q\circ f_B$, and $q$ is the quotient map of Theorem \ref{T:onto}.
 Notice that $\rho$ is a quotient
map, using the equivalence relation defined by $x \not \sim y$  if
there exists uncountable many disjoint cut pairs $\{s,t\} \subset f_B(\bar B)$ separating $f_B(x)$
from $f_B(y)$.  Since the map $f_B$
is canonical up to isotopy, this is an equivalence relation, and so $\rho$ is a $H$ quotient.

We first show that the action of $H$ on this quotient $S^1$ is a
convergence action. Let $(g_i)$ be a sequence of distinct elements
of $G$.  Passing to a subsequence we find $n,p \in \bd X$ such
that for any $x \in \bd X$, $g_i(x) \to p$ and $g_i^{-1}(x) \to
n$. Since $G$ is rank 1, every non-empty open set of $\bd X$ contains points at infinite Tits distance from $n$.  Every arc of $S^1$ contains the image of a non-empty open subset of $\bd X$.  Thus every open arc of $S^1$ contains the image of a point at infinite Tits distance from $n$.

 Let $\hat p = \rho(p)$ and $\hat n = \rho (n)$.  For $v,u \in S^1$ we define the Tits distance
$d_T(u,v) =d_T(\rho^{-1}(u),\rho^{-1}(v))$.  Since the action of $G$ preserves Tits distance on $\bd X$, it also preserves it on $S^1$.  Let  $U$, open arc about of $\hat p$.  There exists  $\epsilon>0$  and, such that $d_T(\hat p,J) >\epsilon$  where $J = S^1-U$.

 Let $K$ be a
closed arc of $S^1$ with endpoints $a$ and $b$ and $\hat n \not \in K$.  It suffices to show
that $g_n(K)\subset U$ for all $i >>0$.  If $\rho^{-1}(K) \cap B_T(n, \pi) = \emptyset$ then by $\pi$-convergence, $g_i(\rho^{-1}(K)) \to \bar B_T(p,\epsilon)$ so $g_i(K) \subset U $ for all $i>>0$.   Thus we may assume that $K \subset \rho( \bar B_T(n,\pi))$.  In fact we can now reduce to the case where $K = \rho(\alpha)$ where $\alpha$ is a Tits arc of length less that $\frac \epsilon 2$.
The arc $\alpha$ contains the image $q$ of a point at infinite Tits distance from $n$.  Thus $g_i(q) \to p$.
By the lower semi-continuity of the Tits metric on $ \bd X$, for $i >>0$, $d_T(g_i(q), J) > \frac \epsilon 2$,
and it follows that $g_i(K) \subset U$.

Thus $H$ acts as a convergence group on $S^1$, so $H$ acts properly
on $\HH^2$ with  the action of $H$ on $S^1$ coming the induced
action of $H$ on $\bd \HH^2 = S^1$.  Clearly for any gap $J$ of
$B$, $\pi(J)$ is a single point.  Thus $stab(J) \cap H<
stab(\pi(J))$.  Since the stabilizer of single boundary point in a
Fuchsian group is either finite or peripheral, the proof is
complete.
\end{proof}

\begin{Cor}\label{C:necstab} If $I$ is an interval of $\cP$ which contains
a necklace $N$, then $stab(I)$ is virtually cyclic. If $N$ is
interior to $I$ then $stab(I)$ is finite.
\end{Cor}
\begin{proof}
In any case the $stab(I)$ will stabilize $N$ and a gap of $N$, and
thus be virtually cyclic. If $N$ is interior to $I$, then the $stab(I)$
stabilizes $N$ and two distinct gaps of $N$, and so $stab(I)$ is
finite.
\end{proof}
\begin{Cor} If $G$ leaves a necklace invariant, then $G$ is virtually a closed hyperbolic surface group.
\end{Cor}
\begin{proof}
By Theorem \ref{T:necklace}, $G$ is virtually a Fuchsian group.  Since $G$ has one end, $G$ is virtually a closed hyperbolic surface group.
\end{proof}

{\bf For the remainder of the paper, we assume (in addition) that $G$ is not virtually a closed surface group, so $G$ acts on $T$ without global fixed point.}

\subsection{T is discrete}

The tree $T$ we have constructed from the boundary of $X$ is
a-priori a possibly non discrete $\Bbb R-$tree. We recall that a
$G$-tree is termed minimal if it does not contain a $G$-invariant
proper subtree.

\begin{Lem} The tree $T$ is minimal.
\end{Lem}
\begin{proof} Suppose not, then we pass to
the minimal subtree contained in $T$, which we denote by $T_m$.
Let $\cP _m=\cP \cap T_m$ (i.e. $\cP _m $ are the elements of $\cP
$ that correspond to points of $T_m$).  Let $A$ be a necklace or
an inseparable cut pair. Then there are open sets $U$ and $V$ of
$\bd X$ separated by a cut pair contained in $A$ with $U \cap A= \emptyset$.  Since the action of $G$ on $X$ is rank
1, there is a rank 1 hyperbolic element $h$ with $h^- \in U$ and
$h^+ \in V$.  It can be shown using $\pi$-convergence  that $h$ acts by
translation on a line  $L \subset T$ with $A \in L$.  Thus $A \in \cP_m$ and so every
cut pair and necklace is in $\cP_m$.

Recall by definition that we removed the terminal points of $\cP$ before we glued in intervals to form $T$.  Notice that any non-terminal maximal inseparable set is between two elements which are either
necklaces or inseparable cut pairs.  Since all necklaces and inseparable cut pairs are in $\cP_m$,
all non-terminal elements of $\cP$ are in $\cP_m \subset T_m$.  It follows that  $T_m =T$.
\end{proof}
\begin{Thm} $T$ is discrete.
\end{Thm}
\begin{proof}   We first show that
$T$ is stable.  Let $I$ be an interval of $\cP$ with
infinitely many elements.  We will show that $stab(I)$ is a finite
group. If $I$ contains a necklace in its interior then $stab(I)$
is finite by Corollary \ref{C:necstab}. Otherwise there are
infinitely many inseparable cut pairs in $I$.

Assume that $stab(I)$ is infinite and let $g_n\in stab(I)$ be an
infinite sequence of distinct elements. By passing to a
subsequence we may assume that $g_n(x)\to p$, and $g_n^{-1}(x)\to n$ for
some (all) $x\in X$. There are inseparable cut pairs $A, B  \in I$
with $A$ separating $B$ from $p \in \partial X$. Choose continua
$Y$, $Z$ with $Y \cup Z = \bd X$ and $Y \cap Z = B$.
 We may assume that $A \subset Y$, which forces $p\in Y$.

 Since $G$ is rank 1, there exists  $z \in Z$ with  $d_T(z, n) = \iy$.
 By $\pi$-convergence $g_n(z) \to p$.
 Consider
 $V = \bigcup\limits_{n\ge 0} h^n(Z)$.  The set $V$ is connected since
each set in the union contains $B$.
 Also $A \not \in V$.  But $p \in \bar V$, the closure of $V$.  This
contradicts the fact that $A$ separates $B$ from $p$.
  Thus $stab(I)$ is  finite.

  Since there is a bound on the size of a finite subgroup in $G$, the action
of $G$ on $T$ is stable.

 Suppose that $T$ is not discrete.  Then by \cite{Le} there is an
$\R$-tree $S$
 and a $G$-invariant  quotient map $f:T \to S$  such that $G$ acts
non-trivially  on $S$ by isometries.
 Furthermore for each non-singleton arc $\alpha \subset S$,
$f^{-1}(\alpha)\cap \cP$ is an infinite interval
 of $\cP$, and the stabilizer of an arc of $S$ is the stabilizer of an
arc of $T$.
 It follows that the arc stabilizers of $S$ are finite and that the action
of $G$ on $S$ is stable.

 Since $G$ is one-ended, by the Rips machine applied to $S$,
 there is a  element $h\in G$ which acts by translation on a line in $S$
such that $G$ virtually splits over $h$.
 This implies that $\{h^\pm\}$ is a cut pair, however since $h$ acts by
translation on a line in $S$,
 it acts by translation on a line in $T$, and so $h^\pm$ will correspond
to the ends of this line (by $\pi$-convergence).
 Thus $\{ h^\pm\}$ is not an element of $\cP$, nor is it a subset of a
necklace of $\cP$.  This contradicts the fact that $\{h^\pm\}$ is
a cut pair.
 Thus $T$ is discrete.
 \end{proof}

 \begin{Lem} \label{L:ehelper} Let  $A\in \cP$ be an inseparable cut pair
and $C\in \cP\cap T$ adjacent to $A$.
  Then $A$ doesn't separate $C$ from a point of $\Lambda  stab[C,A]$.
  Furthermore $\Lambda  stab[C,A]$ is contained in $C$.
 \end{Lem}
 \begin{proof}
 Suppose that $p \in \Lambda  stab[C,A]$ is separated from $C$ by $A$.
Choose a $(g_i) \subset stab[C,A]$ with $g_i(x) \to p$ and
$g_i^{-1}(x)\to n$ for some $n \in \bd X$.
 Since $C\in \cP \cap T$, there is a necklace or cut pair $E\in \cP$ with $C
\in (E, A)$.  Since $G$ is rank one, there is a point $q$ at
infinite Tits distance from $n$ separated from $C$ by $E$.   It
follows that $g_i(q) \to p$ but this contradicts $A$ separating
$C$ from $p$.

To show that $\Lambda  stab[C,A]$ is contained in $C$ we argue
similarly. Suppose $p \in \Lambda  stab[C,A]$ does not lie in $C$.
Choose a $(g_i) \subset stab[C,A]$ with $g_i(x) \to p$ and
$g_i^{-1}(x)\to n$ for some $n \in \bd X$. Since $A$ is not a
terminal point of $T$, there is a point $q$ at infinite Tits
distance from $n$ separated from $p$ by $A$. It follows that
$g_i(q) \to p$. We note however that $p$ is separated from $C$ by
a cut pair $c,d$. So we can write $\bd X=E\cup F$ with $E,F$
continua such that $E\cap F=\{c,d\}$. Let's say $A\subset F$,
$p\in E$. Now since $g_i$ fixes $A$, $g_i(q)\in F$ for all $i$. So
$g_i(q)$ does not converge to $p$, a contradiction.

 \end{proof}
 \begin{Thm}\label{T:cutpairst} If $A $ is an inseparable cut pair then
there is a hyperbolic $g$ with $A = \{g^\pm\}$.
 \end{Thm}
 \begin{proof} Suppose not.  Let $C,D$ be elements of $\cP \cap T$ adjacent to
$A$. Since $G$ does not split over a
finite group, there are infinite finitely generated subgroups  $H=
\langle h_1\dots h_n\rangle$ and $\hat H = \langle \hat h_1\dots
\hat h_m\rangle$ stabilizing $\{C,A\}$ and  $\{A, D\}$
respectively.

By passing to subgroups of index 2, if necessary,
we may assume that $H$ and $\hat H$ fix $A$ pointwise.

It follows from \cite{SWE} that $A \subset  \Lambda Z_{h_i} $ and
$A \subset  \Lambda Z_{\hat h_j}$ for all $i,j$. Since $Z_{h_i}$ and
$Z_{\hat h_j}$ are convex, it follows from \cite[Theorem 16]{SWE}
that $Z_H = \bigcap\limits_i Z_{h_i}$ and $Z_{\hat H}
=\bigcap\limits_i Z_{h_i}$ are convex and that $A \subset   \Lambda
Z_H$ and $A \subset  \Lambda Z_{\hat H}$.  Thus by \cite[Theorem
16]{SWE} $Z = Z_H \cap Z_{\hat H}$ is convex and $A \subset  \Lambda
Z = \Lambda Z_H \cap \Lambda Z_{\hat H}$.

Since $Z$ is convex, it is a CAT(0) group, and by \cite[Theorem
11]{SWE} there is an element $g \in Z$ of infinite order
\cite[Corollary to Theorem 17]{SWE}.  Since $H, \hat H$ are infinite,
$\Lambda H, \Lambda \hat H$ are non-empty.

 By Lemma \ref{L:ehelper},
$\Lambda H $ is not separated from $C$ by $A$ and similarly $\Lambda \hat H $ is not separated from $D$ by $A$.  By Lemma
\ref{L:centlim}, $g$ fixes $\Lambda H$ and $\Lambda \hat H$. If either $\Lambda H \subset A$ or
$\Lambda H \subset A$, we are done by Lemma \ref{L:limsets}.
 Either $\Lambda H \subset C$, in which case $g(C) =C$, or  points of $\Lambda H$ are separated from $D$ by $C$, and similarly for $\Lambda \hat H$ and $D$.

  It could a priori happen that $g$ acts by translation on a line of $T$ which contains $A,C,D$.  In that case however, $g$ is rank 1 (see proof of Corollary \ref{C:rank11}).  By Lemma \ref{L:limsets} $|\Lambda H|, |\Lambda \hat H|\ge 2$ Thus $\Lambda H \cup \Lambda \hat H \not \subset \{g^\pm\}$ contradicting the fact that $g$ fixes both and is rank 1.

The only remaining possibility is that  $g(C) =C$ and  $g(D) =D$,
so $g(A) =A$.     It follows by Lemma \ref{L:ehelper} that
$g^{\pm}\subset C$ and $g^{\pm}\subset D$. Therefore $g^{\pm}\in
C \cap D = A$ and so $\{g^\pm\} =A$.
   \end{proof}

\begin{Cor} Let $A,B,C \in \cP \cap T$ with $B \in (A,C)$ and $B$ inseparable
cut pair. If $A$ and $C$ are both adjacent to $B$ then the
stabilizer of the interval $[A,C]$ is finite or virtually $\langle
g\rangle$ where $g$ is hyperbolic with $\{g^\pm \} = B$.
\end{Cor}
\begin{proof}  Let $H = stab[A,C]$, so $H=stab[A,B] \cap
stab[B,C]$. By Lemma \ref{L:ehelper}, $\Lambda H \subset A$ and
$\Lambda H \subset C$, so $\Lambda H \subset A \cap C = B$.  If
$H$ is not finite, it follows that $H$ is  virtually $\langle
g\rangle$ where $g$ is hyperbolic with $\{g^\pm \} = B$.
\end{proof}
Let $I \subset \cP \cap T$ be an interval of $\cP$.  Notice that if
$I$ contains three elements of $\cP$ which are not inseparable cut
pairs, then $stab(I)$ is finite.

\subsection{$G$ splits over a 2-ended group}

We recall here some terminology for group actions on trees and
graphs of groups.

If $\G $ is a graph of groups and $V$ is a vertex group of $\G $
then we say that we can refine $\G $ at $V $ if the following
holds: There is a graph of groups decomposition of $V$ such that
all edge groups of $\G $ adjacent to $V$ are contained in vertex
groups of the graph of groups decomposition of $V$. We obtain the
refinement of $\G $ by substituting $V$ by its graph decomposition
and joining the edges of $\G $ adjacent to $V$ to the vertex
groups of the graph containing the corresponding edge groups. A
vertex $v$ of a graph of groups labeled by a group $V$ is called
non-reduced if it is adjacent to at most 2 edges, none of which is
a loop, and all edges adjacent to $v$ are labeled by $V$. A graph
of groups which contains no such vertices is called reduced. We
say that a subgroup $H<G$ is elliptic in $\G $ if it is contained
in a conjugate of a vertex group of $\G $.

We proceed now to show that $G$ splits over a 2-ended group. We
will show this considering the action of $G$ on $T$. In the next
section we show that this action gives also the JSJ decomposition
of $G$. We remark that $T$ has at least one vertex that
corresponds either to an inseparable cut pair or to a necklace. We
will show in each case that one can get a splitting of $G$ over a
2-ended group.

 We consider the graph of groups $\G $ obtained from the action
of $G$ on $T$.

By Theorem \ref{T:necklace} if a vertex $v$ of $\G $ corresponds
to a necklace then the vertex is labeled either by virtually a
subgroup of $\Bbb Z ^2$ or it is a fuchsian group. In the second
case, $\G $ has an edge labeled by a peripheral subgroup of the
fuchsian group. In the first case, if the vertex is labeled by
$\Bbb Z ^2$, then the limit set of the vertex stabilizer is a
circle. Since this limit set is a necklace, we have that the
boundary of $G$ is a circle so $G$ is virtually $\Bbb Z ^2$.
Otherwise the vertex is labeled by virtually $\Bbb Z$. If the
vertex corresponds to a branch point, we see that all edges
adjacent to it are labeled by 2-ended groups and the graph is
reduced at this vertex. Otherwise a vertex $u$ adjacent to it is
labeled by a group containing the stabilizer of a maximal
inseparable set which is a branch point of $T$. The edge with
endpoints $[u,v]$ is labeled by a 2-ended group and gives a
splitting of $G$.

Let's now consider a vertex $v$ of $T$ corresponding to an
inseparable pair $\{a,b \}$. We consider first the case that the
stabilizer of $v$ is a two-ended group. If $v$ corresponds to a
branch point of $T$, then, since the edges adjacent to $v$ are
labeled by subgroups of the stabilizer of $v$, all these edges give
reduced splittings of $G$ over two ended groups. If $v$ does not
correspond to a branch point, then an edge adjacent to $v$
corresponds to a branch point of $T$. This edge is labeled by a
two ended group and gives a reduced splitting of $G$ over a two
ended group.

Let $V$ be the stabilizer of the vertex $v$ corresponding to
$a,b$. By theorem \ref{T:cutpairst} there is a hyperbolic element
$g$ in $V$ with $\{g^{\infty }, g^{-\infty} \}=\{a,b\}$. By the
algebraic torus theorem \cite {D-S} to show that $G$ splits over a
2-ended group it suffices to show that $X/<g>$ has more than one
end.

Let $\bar X=X\cup \partial X$. Let $L$ be an axis for $g$ and let
$D$ be the translation length of $g$. Let $U,V$ be small
neighborhoods around $a,b$ respectively. To fix ideas we define
$U,V$ as follows: we take a base point $x$ on $L$ and we consider
the ray from $x$ to $a$. Let's call this ray $c_a$. Let's call
$c_y$ the (possibly finite) ray from $x$ to a point $y\in \bar X$.
We define $U$ as follows:
$$ U=\{y\in \partial X: d(c_y(10D),c_a(10D))<1 \}$$
$V$ is defined similarly. We claim that there is a $K$ such that a
$K$-neighborhood of $L$ union $U \cup V$ separates any two
components of $\partial X-{a,b}$.  Suppose that there are
components $C_1,C_2$ of $\partial X-\{a,b\}$ such that for any $n$
there is a component $T_n$ of $\bar X$ minus $U\cup V \cup N_n(L)$
such that $T_n$ intersects both $C_1,C_2$. Without loss of
generality we can assume that the $T_n's$ are nested. Then if $T
=\bigcup T_n$, $T\subset \bd X$ is closed, connected and
intersects both $C_1,C_2$. This is a contradiction since $T$ does
not contain $a$ or $b$.

We recall the fact (see \cite{ONT}) that geodesics in $X$ are
`almost extendable' i.e. there is an $A>0$ such that if $[a,b]$ is
a finite geodesic in $X$ then there is an infinite ray $[a,c]$
($c\in
\partial X$) such that $d(b,[a,c])\leq A$.

Let $p:X\to L$ be the projection map from $X$ to $L$.

We claim that there is an $M>0$ such that the $M$-neighborhood of
$L$ separates $X$ in at least two components $Y_1,Y_2$ such that
$Y_i$ is not contained in any finite neighborhood of $L$
$(i=1,2)$. This claim implies that $G$ either splits over a
2-ended group or is virtually a surface group by \cite{D-S}.

We now prove the claim. Let $X_1,X_2$ be two distinct components
of $\partial X- \{a,b\}$. We claim that there are two infinite
rays $r_1,r_2$ from $x$ to $X_1,X_2$ respectively which are
perpendicular to $L$. Indeed we consider rays from $x$ to ,say,
points on $\bar X_1$. The angle that these rays form with $L$
varies continuously and takes the values $\{0,\pi \}$ at
$\{a,b\}$. Therefore some such ray $r_1$ is perpendicular to $L$.
We argue similarly for $r_2$.

Given $n>A, n\in \mathbb{N}$ let $R_1=r_1(n+1),R_2=r_2(n+1)$. If
$R_1,R_2$ are not contained in the same component of $X-N_n(L)$
then the claim is proven. Otherwise there is a path $p$ in
$X-N_n(L)$ joining $R_1$ to $R_2$. For every $y\in X$ we consider
$p(y)\in L$ and we pick an infinite ray $r_y$ from $p(y)$ such
that $d(r_y,y)\leq A$. For $R_1,R_2$ we choose the corresponding
rays to be $r_1,r_2$ respectively. Clearly there are $y_1,y_2\in
p$ such that $d(y_1,y_2)<1$ and the corresponding rays
$r_{y_1},r_{y_2}$ define points in distinct components of
$\partial X- \{ a,b \} $. Since $<g>$ acts cocompactly on $L$ we
may translate $p(y_1),p(y_2)$ close to $x$, say
$$d(g^k(p(y_1)),x)<2D, d(g^k(p(y_2)),x)<2D$$
Then $g^k(r_{y_1}),g^k(r_{y_2})$ define points in distinct
components of $\partial X- \{ a,b \} $ and if $n>K$ these two
points are not separated by the $K$-neighborhood of $L$ union $U
\cup V$. This is a contradiction.

We have shown the following:

\begin{Thm}\label{T:split} Let $G$ be a one ended group acting
geometrically on a CAT(0) space $X$. If a pair of points $\{a,b
\}$ separates $\bd X$ then either $G$ splits over a 2-ended group
or $G$ is virtually a surface group.
 \end{Thm}

\subsection{JSJ decompositions}

\begin{Def}
Let $G$ be a group acting on a tree $T$ and let $\G $ be the
quotient graph of groups. We say that this action (or the
decomposition $\G $) is canonical if for every automorphism
$\alpha $ of $G$ there is an automorphism $C_{\alpha }$ of $T$
such that $\alpha (g)C_{\alpha }=C_{\alpha }g$ for every $g\in G$.
\end{Def}
Bowditch (\cite{BOW}) constructed a canonical JSJ-decomposition
for hyperbolic groups. Swarup-Scott (\cite {S-S}) constructed a
canonical JSJ-decomposition for finitely presented groups in
general. We show here how to deduce a canonical JSJ-decomposition
for CAT(0) groups using their CAT(0) boundary. Our
JSJ-decomposition is similar to the one in \cite{S-S}.

To describe our decomposition we use the notation of the previous
sections. So $G$ is a CAT(0) group acting on a CAT(0) space $X$.
From the pairs of cut points of the boundary of $X$, we construct
an $\Bbb R$-tree $T$ on which $G$ acts non-trivially (unless $G$
has no splittings over two-ended groups or $G$ is virtually of the
form $H\times \Bbb Z$, in these cases our JSJ decomposition is
trivial). As we showed in the previous section, $T$ is discrete.
Then the graph of groups decomposition $\G $ given by the quotient
graph of the action of $G$ on $T$ is the JSJ decomposition of $G$.

Vertex groups of this graph are either fuchsian or groups that
have no splitting over 2-ended groups in which the adjacent edge
groups are elliptic or groups that stabilize inseparable cut
pairs. Edges incident to fuchsian groups are labeled by 2-ended
groups by Theorem \ref{T:necklace}. The only other type of edges
are edges joining an element $C$ of $\cP$ corresponding to an
inseparable subset of $\bd X$ and an inseparable pair $\{a,b \}$
of $\cP$. These edges are not labeled necessarily by 2-ended
groups. So our decomposition might differ from the one obtained by
Rips-Sela \cite{R-S} or Dunwoody-Sageev \cite {D-Sa}. We give an
example to illustrate this point. Let
$$A=\Bbb F_4\times \Bbb Z=<a,b,c,d>\times <x>$$
$$B=\Bbb F_2\times \Bbb F_2=<a,b>\times <x,y>$$
$$C=\Bbb F_2\times \Bbb F_2=<c,d>\times <x,z>$$
and consider the group given by the graph of groups
$$G=B*_EA*_DC$$ where
$E=<a,b>\times <x>=\Bbb F_2\times \Bbb Z$, $D=<c,d>\times <x>=\Bbb
F_2\times \Bbb Z$. Then the JSJ decomposition that we obtain is
the graph of groups $B*_EA*_DC$. However we remark that this
decomposition can be refined to give the ordinary JSJ
decomposition by splitting $A$ over $x$, i.e. we have
$$G=B*_EA*_DC=B*_E*(E*_{<x>}D)*_DC=A*_{<x>}C$$
where the latter decomposition is the ordinary JSJ decomposition.

We claim that one can always refine our decomposition to a JSJ
decomposition over 2-ended groups given in \cite {D-Sa}. Let $T_J$
be the tree corresponding to the JSJ decomposition in \cite
{D-Sa}. To show that $\G $ can be refined, it is enough to show
that edge groups of $\G $ act elliptically on $T_J$. This is
clearly true for the 2-ended vertex groups. Now let
 $V$ be the stabilizer of an inseparable cut pair $\{a,b
\}$ and let $e$ be an edge joining $\{a,b \}$ to an element $C$ of
$\cP$ corresponding to an inseparable subset of $\bd X$. If $E$ is
the edge stabilizer of $e$ by lemma \ref{L:ehelper} $\Lambda E$ is
contained in $C$. If $E$ acts hyperbolically on $T_J$, then there
is a hyperbolic element $h\in E$ and the 2 points of $\Lambda
\langle h \rangle$ are separated by any pair of points of $\bd X$
corresponding to limit points of an edge on the axis of $h$ on
$T_J$. This is a contradiction since $\Lambda \langle h\rangle
\subset C$. It follows that any edge group of $\G $ fixes a point
of $T_J$. We remark that the only vertex groups of $\G $ that
might act hyperbolically on $T_J$ are vertex groups that stabilize
inseparable cut pairs. To refine $\G $ we let all such groups act
on $T_J$ and we replace them in $\G $ by the decompositions we
obtain in this way. Doing this refinement and collapsing we obtain
a JSJ-decomposition of $G$ over 2-ended groups in the sense of
Dunwoody-Sageev.

In order to show that our JSJ decomposition is canonical, we
observe that if $G$ splits over a 2-ended subgroup $Z$, then the two
limit points of $Z$ form  a cut pair of $\partial X$. So if
$\alpha $ is an automorphism of $G$ the limit points of $\alpha
(Z)$ is also a separating pair for $\bd X$. This observation
suffices to show that our JSJ decomposition is canonical. We now
explain this in detail.

Let $\alpha $ be an automorphism of $G$. If $H$ is a fuchsian
hanging subgroup of $\G $, then the limit set of $H$ in $\bd X$ is
a necklace and $H$ fixes the vertex of $T$ corresponding to this
necklace. $\alpha (H)$ is also a hanging fuchsian group whose
limit set is a necklace and $\alpha (H)$ fixes the corresponding
vertex. A vertex $V$ of $\G $ which is virtually of the form
$H\times \Bbb Z$  is the stabilizer of an inseparable pair $\{a,b
\}$. Clearly $\alpha (V)$ is also a stabilizer of an inseparable
pair so it fixes also a vertex of $T$. Finally we remark that
$\pi $ convergence implies that the stabilizer $V$ of a maximal
inseparable set has a limit set which is contained in the maximal
inseparable set. So $\alpha (V)$ fixes also a vertex of $T$.

We observe that adjacent vertex groups are mapped to adjacent
vertex groups by $\alpha $. This is because splittings $A*_CB$ are
mapped to splittings $\alpha (A)*_{\alpha (C)}\alpha (B)$ by
$\alpha $.

We have shown the following:
\begin{Thm}\label{T:jsj} Let $G$ be a one ended group acting
geometrically on a CAT(0) space $X$. Then the JSJ-tree of $\bd X$
is a simplicial tree $T$ and the graph of groups $T/G$ gives a
canonical JSJ decomposition of $G$ over 2-ended groups.
\end{Thm}

\end{document}
\bye